\documentclass[12pt]{article}
\input{epsf}
\usepackage{amsmath}
\usepackage{amsfonts}
\usepackage{amscd}
\newtheorem{definition}{Definition}[section]

\newtheorem{theorem}[definition]{Theorem}

\newtheorem{conjecture}[definition]{Conjecture}
\def\block{\hbox{${\vcenter{\vbox{\hrule height 0.4pt\hbox{\vrule 
width 0.4pt 
height 6pt \kern 5pt\vrule width 0.4pt}\hrule height 0.4pt }}}$}}
\def\qed{\hfill\block\medskip}
\def\deg{\hbox{deg}}
\def\Hilb{\hbox{Hilb}}

\def\Vect{\hbox{Vect}}
\def\Hom{\hbox{Hom}}
\def\End{\hbox{End}}
\def\N{\mathbb{N}}

\def\To{\Rightarrow}

\def\ol{\overline}
\begin{document}

\title{Finite groups, spherical 2-categories, and 4-manifold invariants}
\author{Marco Mackaay \\ \\ \'Area Departamental de Matem\'atica,\ UCEH, \\ Universidade do Algarve, \\ 8000 Faro, \\ Portugal, \\ e-mail mmackaay@ualg.pt}
\date{March 14, 1999}
\maketitle
\begin{abstract}
In this paper we define a class of state-sum invariants of closed 
oriented piece-wise linear 4-manifolds using finite groups. The definition 
of these state-sums follows from the general abstract construction of 
4-manifold invariants using spherical 2-categories, as we defined in 
an earlier paper. We show that the state-sum invariants of Birmingham and 
Rakowski, who studied Dijkgraaf-Witten type 
invariants in dimension 4, are special examples of the 
general construction that we present in this paper. They showed 
that their invariants are non-trivial by some explicit computations, so 
our construction includes interesting examples already. Finally, we 
indicate how our construction is related to homotopy 3-types. This connection 
suggests that there are many more interesting examples of our construction 
to be found in the work on homotopy 3-types, by Brown, for 
example.  
\end{abstract}

\section{Introduction}

In \cite{Ma99} we defined spherical 2-categories and showed how to construct state-sum invariants 
of closed oriented PL 4-manifolds with them. Roughly speaking spherical 
2-categories are monoidal 2-categories with 
duals, as defined by Baez and Langford~\cite{BL982,BL98}, such that the 
categorical trace satisfies a small set of conditions. The main point in that 
paper was to find a construction that would generalize Crane and Frenkel's construction~\cite{CF94}, which uses involutory Hopf categories, and Crane and Yetter's construction ~\cite{CKY97, CY93}, which uses tortile categories. Let us 
exlain this in some detail.

\begin{figure}
\begin{picture}(380,200)
\put(380,0){\makebox(0,0){2D}}
\put(30,0){\makebox(0,0){algebra}}
\put(30,10){\line(0,1){30}}
\put(380,50){\makebox(0,0){3D}}
\put(30,50){\makebox(0,0){Hopf algebra}}
\put(30,60){\line(0,1){30}}
\put(380,100){\makebox(0,0){4D}}
\put(30,100){\makebox(0,0){trialgebra}}
\put(70,50){\line(1,0){30}}
\put(70,100){\line(1,0){30}}
\put(155,100){\makebox(0,0){Hopf category}}
\put(155,50){\makebox(0,0){monoidal category}}
\put(200,100){\line(1,0){30}}
\put(295,100){\makebox(0,0){monoidal 2-category}}
\put(155,60){\line(0,1){30}}
\put(70,10){\line(2,1){30}}
\put(200,70){\line(2,1){30}}
\end{picture}
\caption{The categorical ladder}
\label{catlad}
\end{figure}
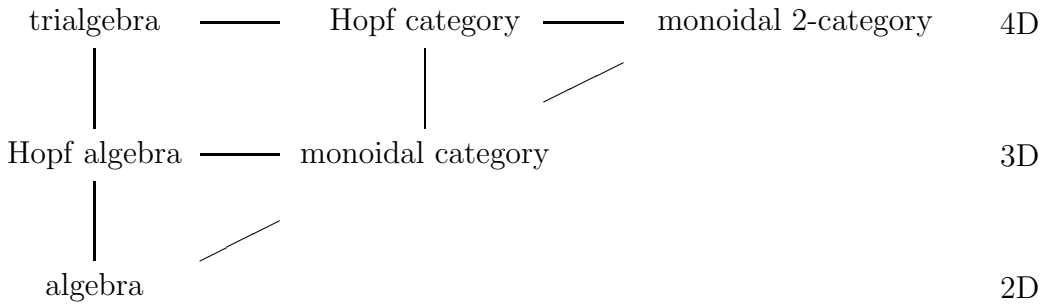

\begin{figure}
\centerline{
\epsfbox{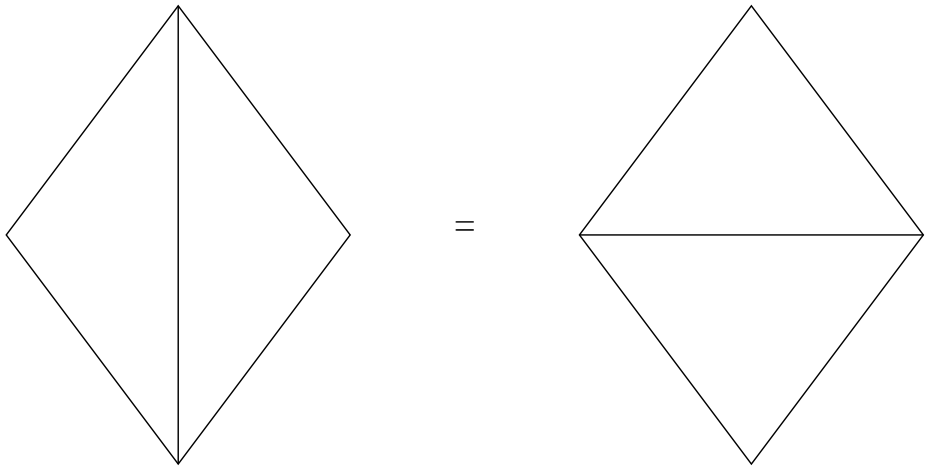}
}
\caption{$2 \rightleftharpoons 2$ Pachner move}
\label{pach22}
\end{figure}

\begin{figure}
\centerline{
\epsfbox{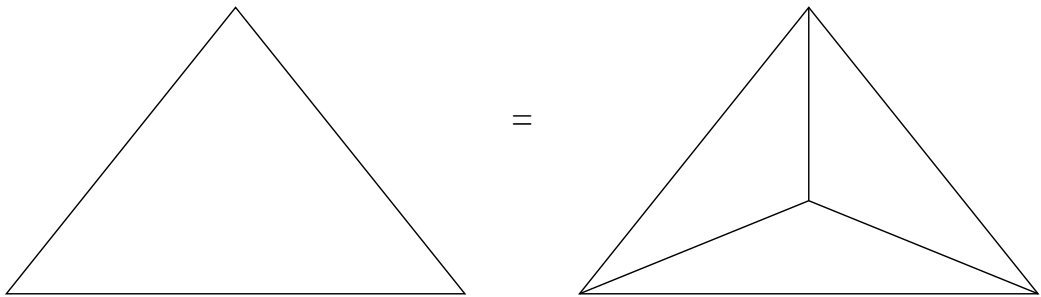}
}
\caption{$1 \rightleftharpoons 3$ Pachner move}
\label{pach13}
\end{figure}

In \cite{CF94} Crane and Frenkel sketched a 
general approach to the construction of 4D TQFT's. At its core is the 
so-called {\it categorical ladder}, as shown in Fig.~\ref{catlad}. 
For an explanation of this diagram see~\cite{CF94}. 
Here we only explain the diagonal going from {\it algebra} to {\it monoidal 
2-category}. In \cite{FHK94} the authors show that any finite dimensional 
semi-simple associative algebra, $A$, can be used for the construction of 
state-sum invariants of 2-manifolds (surfaces). This is not the place to 
recall the construction in detail. What is of interest to us in this 
introduction is the main idea behind the construction, which accounts for 
the invariance of the state-sums. To show that the value of a state-sum of a 
triangulated surface does not depend on the chosen triangulation, so that 
this value is a topological invariant of the surface, one has to prove that 
the value does not change under any of the 2D Pachner moves. 
The $n$-dimensional Pachner moves express the combinatorial 
equivalence relation between triangulations of two PL homeomorphic closed 
oriented PL $n$-manifolds. More precisely, two triangulated $n$-manifolds, 
$(M_1,{\cal T}_1)$ and $(M_2,{\cal T}_2)$, are PL homeomorphic if and only if 
there exists a finite sequence of $n$-dimensional Pachner moves which 
transform ${\cal T}_1$ into a new, but PL-homeomorphic, triangulation of 
$M_1$ which is isomorphic to ${\cal T}_2$ as a simplicial complex. In 
dimension 2 the Pachner moves are the ones depicted in Fig.~\ref{pach22} and 
Fig.~\ref{pach13}.

\begin{figure}
\centerline{
\epsfbox{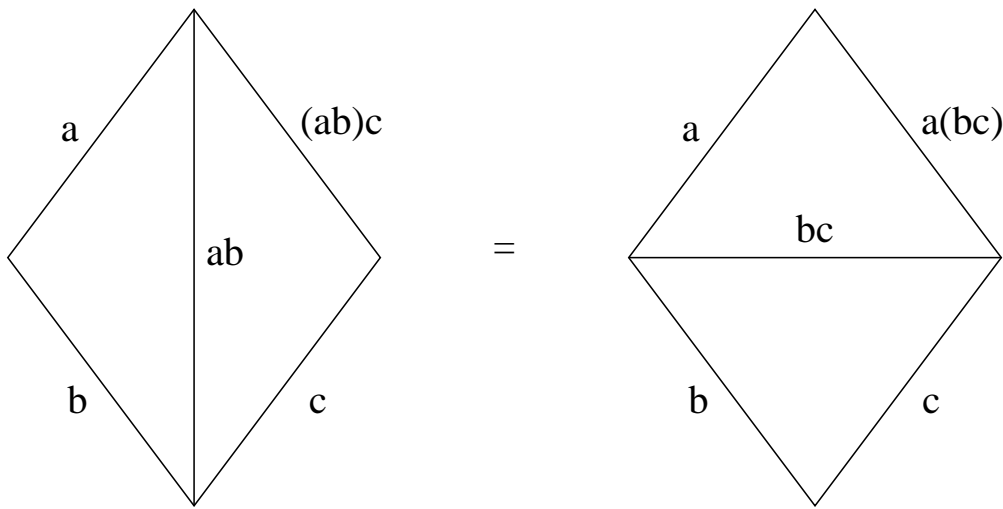}
}
\caption{$2 \rightleftharpoons 2$ Pachner move}
\label{pach222}
\end{figure}

As shown in Fig.~\ref{pach222} the $2\rightleftharpoons 2$ move can be 
interpreted as a diagrammatic way of expressing the associativity of $A$, if 
one labels the edges with elements of $A$. 
For the partition function of the 
state-sum of a triangulated surface $(S,{\cal T})$ one has to label the edges 
of $\cal T$ with basis elements of $A$ and associate to each triangle the 
multiplication constant of $A$ determined by the three labels on the edges 
in the boundary of the triangle. Therefore, Fig.~\ref{pach222} shows that invariance 
of the state-sum under the $2\rightleftharpoons 2$ Pachner move follows from 
the equations satisfied by the multiplication constants expressing the  
associativity of $A$. Invariance under the other 2D move is a bit trickier 
and involves the semi-simplicity of $A$ as well. Since we are only interested 
in the foundational ideas in this introduction, we do not explain the 
invariance under this move.

In dimension 3 there are two different constructions of state-sum invariants. 
One is due to Kuperberg~\cite{Ku91}, who uses finite 
dimensional involutory Hopf algebras as algebraic input. In~\cite{CFS94} the 
reader can find a detailed account of Kuperberg's invariants, which we do 
not reproduce. The other construction in dimension 3 is due to Turaev and 
Viro~\cite{TV92}. Several authors~\cite{BW96,Tu94,Ye94} 
generalized Turaev and Viro's original construction, which uses a 
non-degenerate quotient of the 
monoidal category of finite dimensional representations of $U_q(sl(2))$, 
the quantum group corresponding to the Lie algebra $sl(2)$, for certain 
roots of unity $q$, and axiomatized the categorical input that is needed. 
The most ``economical'' axiomatization is due to Barrett and Westbury. They 
defined the notion of a {\it spherical category}, which is 
a monoidal category with duals satisfying some extra conditions. The relation 
with Kuperberg's work was made explicit by Barrett and Westbury~\cite{BW95}: 
the monoidal category of finite dimensional representations of a finite 
dimensional involutory Hopf algebra, $H$, is a particular kind of spherical 
category, $C_H$, and Kuperberg's invariants using $H$ are equal to 
Barrett and Westbury's invariants using $C_H$. However, there are many 
spherical categories whose objects are not representations of 
involutory Hopf algebras, and some of these give very interesting invariants, 
such as the Turaev-Viro invariants. 

\begin{figure}
\centerline{
\epsfbox{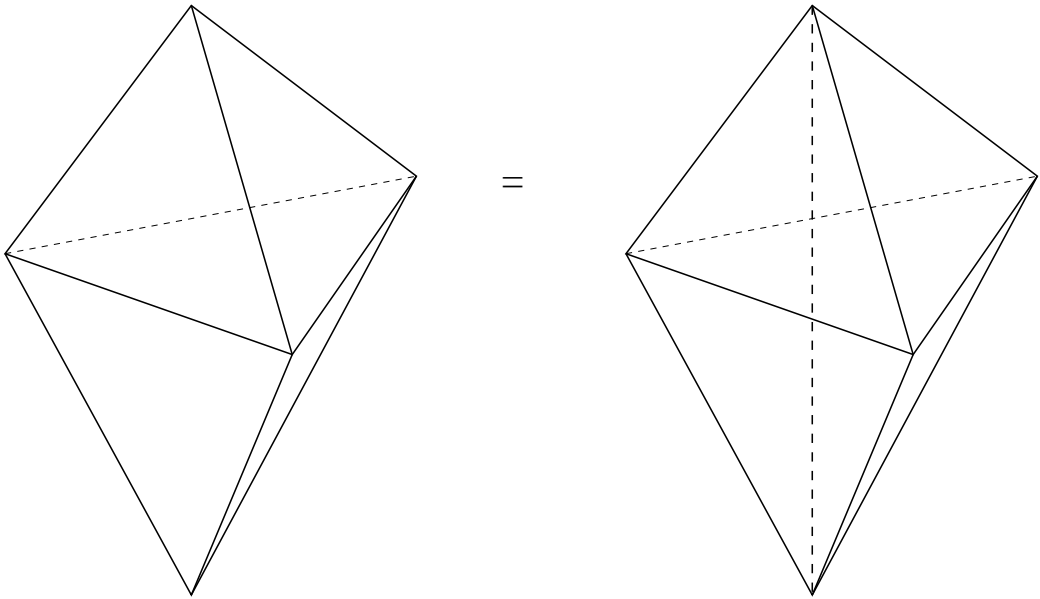}
}
\caption{$2 \rightleftharpoons 3$ Pachner move}
\label{pach23}
\end{figure}

\begin{figure}
\centerline{
\epsfbox{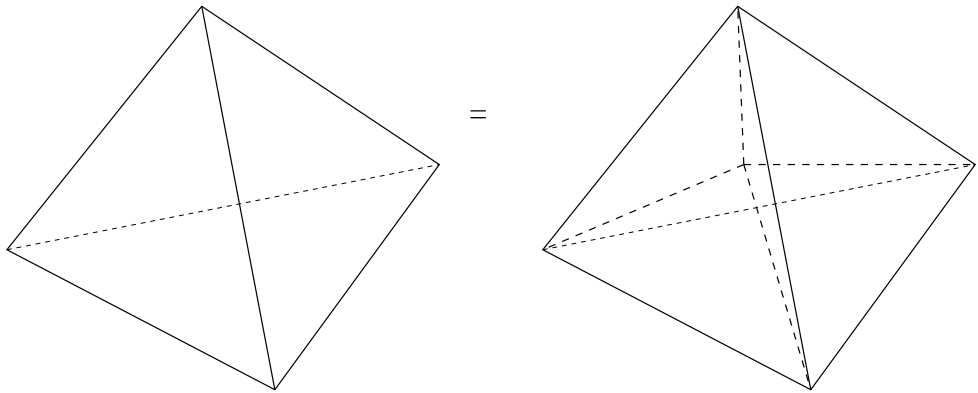}
}
\caption{$1 \rightleftharpoons 4$ Pachner move}
\label{pach14}
\end{figure}

Without 
recalling the details of Barrett and Westbury's construction, let us explain 
how 
this construction of 3D state-sums can be understood in the light of the 
former construction of 2D state-sums. We already remarked that the invariance 
of the 2D state-sums is mainly due to the deep correspondence between the 
combinatorics of the 2D Pachner moves and the algebraic equation 
corresponding to the associativity of algebras. In 3D we have replaced the 
algebra by a monoidal category, which has one extra layer of structure formed 
by the morphisms. The associativity in the algebra, which is an equation, is 
therefore replaced by the so called {\it associator} in the monoidal 
category, which is a natural isomorphism instead of an equation. It is well 
known that the associator of a 
monoidal category has to satisfy a {\it coherence relation}, corresponding 
to the Stasheff pentagon diagram: 
$$
\begin{CD}
X(Y(WZ))@>>> (XY)(WZ)@>>> 
((XY)W)Z\\
@VVV      &           &       @AAA\\
X((YW)Z)  &    @>>>   &      (X(YW))Z
\end{CD}.
$$
In this diagram we have written $XY$ as a shorthand for $X\otimes Y$, and 
the arrows indicate the use of the associator. For the pentagon 
diagram to be commutative the composite of the two instances 
of the associator over the top of the diagram has to be equal to the 
composite of the three 
instances of the associator around the bottom.
Thus, in going from 2D 
to 3D, the elements of an algebra are replaced by the objects in a 
monoidal category and the associativity equation is replaced by the 
associativity isomorphism which 
satisfies a new equation of its own. This ``replacement process'' was called 
{\it categorification} by Crane and Frenkel~\cite{CF94}. There Remains the 
question, what 
the pentagon equation of the associator has to do with the invariance of 
the 3D state-sums. The most satisfactory answer to this question, at least 
in the opinion of the author of this article, is obtained by a close analysis 
of the 3D Pachner moves. These moves can be seen in Figs.~\ref{pach23}, 
\ref{pach14}. 

\begin{figure}
\centerline{\epsfxsize=13cm
\epsfbox{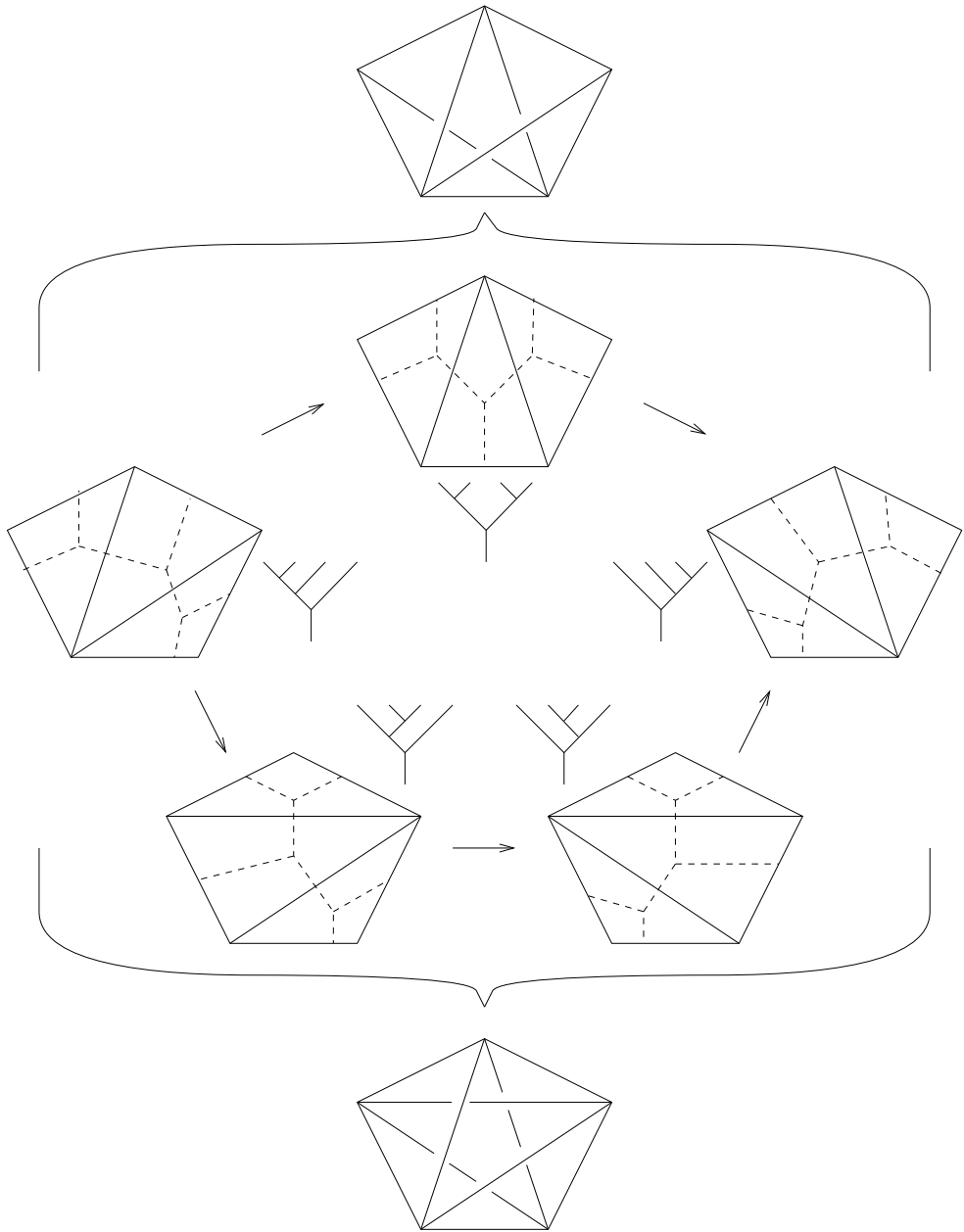}
}
\caption{The pentagon and the $2\rightleftharpoons 3$ Pachner move}
\label{penta}
\end{figure}

In order to understand the invariance of the state-sums 
one would like to see the diagrammatic analogue of categorification and 
its relation with the Stasheff pentagon. This 
diagrammatic interpretation of categorification is due to Carter, Kauffman and 
Saito~\cite{CKS97}. A little thought shows what it should look like: the 
``diagrammatic equations'' expressing the 2D Pachner moves should be 
considered as 
``diagrammatic isomorphisms'', whatever that may be, and a 3D Pachner move 
should be interpreted as a diagrammatic equation between two finite 
sequences of these isomorphisms. It turns out that a 
diagrammatic isomorphism should 
be understood as the gluing of the source and target of the isomorphism. The 
reason for this is that any $n$-dimensional Pachner move corresponds to 
a partition of the boundary of an $n+1$-simplex into two connected 
parts which share a common boundary. Thus, replacing the 
diagrammatic equation that 
corresponds to an $n$-dimensional Pachner move by a 
diagrammatic isomorphism 
can be interpreted as gluing the two $n$-dimensional simplicial complexes,  
which define the two sides of the move, along their common boundary and 
filling up the missing $n+1$-cell in order to obtain the whole $n+1$-simplex. 
This $n+1$-simplex can never be part of the triangulation of an $n$-manifold, 
which is why the $n$D Pachner moves have to be equations in dimension $n$, 
whereas in the triangulation of an $n+1$-manifold there is enough space, 
so the $n$D Pachner moves can no longer hold as equations in dimension $n+1$. 
In order to illustrate this, we have copied Fig.~\ref{penta} from 
\cite{CKS97}, 
with permission from the authors, which 
shows how the $2\rightleftharpoons 3$ Pachner move in 3D can be seen as an 
equation between two finite sequences of 2D Pachner move. 

The arrows indicate 
the gluings which represent the diagrammatic isomorphisms, each of which 
corresponds to a 2D Pachner move, i.e., a tetrahedron. Notice the 
diagrammatic similarity with the Stasheff pentagon! The deep reason why the 
values of the Turaev-Viro type state-sums of 3-manifolds are independent of 
the choice of 
triangulation can now be expressed by saying that the algebraic 
categorification, which is obtained by introducing an associator which 
satisfies the pentagon equation, and the diagrammatic categorification, 
which is obtained by introducing diagrammatic isomorphisms corresponding 
to the 2D Pachner moves which have 
to satisfy the equations corresponding to the 3D Pachner moves, are somehow 
equivalent. Of course this remark is rather vague and the author of this 
article does not know how to make it into a mathematically rigorous 
statement. The problem is to understand the exact relation between the 
coherence relations in (weak) $n$-categories and the $n$-dimensional Pachner 
moves. Unfortunately there are several definitions of weak $n$-categories 
by now, and no one knows whether they are ``equivalent'' in some 
sense. One difference between the various approaches lies in the shape of 
the diagrams that represent the $k$-morphisms for $2\leq k\leq n$; see 
\cite{Ba972} for a nice review of the different approaches. The work 
of Tamsamani~\cite{Ta95}, who defines weak $n$-categories 
via a simplicial approach, might shed some light on the relation between 
coherence relations and Pachner moves one day. Although not mathematically 
rigorous, we hope that the arguments sketched above convince the reader 
that the invariance of Barrett and Westbury's state-sums using spherical 
categories is no miracle. The same arguments also indicate how to proceed 
in dimension 4.           

The deep insight in Crane and Frenkel's paper is that the categorification 
of the 2D state-sums yields the 3D state-sums, and that the invariance of the 
latter are a consequence of the invariance of the former and the general 
principle of categorification. Since they were interested in 
the construction of 4D state-sums, they were led to study the 
categorification of the 3D state-sums. By the arguments above, however vague 
they may seem, the invariance of the 3D state-sums should guarantee the 
invariance of their categorifications. Crane and Frenkel 
chose to categorify Kuperberg's construction, which led them to the 
definition of an involutory Hopf category. This is a monoidal category with 
a comonoidal structure satisfying the axioms of a Hopf algebra up to natural 
isomorphisms, which satisfy a new set of equations themselves. The problem 
is that their categorification inherited Kuperberg's severe restriction of 
the Hopf algebra having to be involutory. As is well known, the most 
interesting 3D invariants are related to the quantum groups, which are not 
involutory. In their conclusions Crane and Frenkel conjecture the possibility 
of categorifying the Turaev-Viro type constructions, which would lead to a 
more general construction, just as Barrett and Westbury's construction is 
more general than Kuperberg's. Crane and Frenkel mention that the 
representations of a Hopf category are categories with a categorified 
module structure, as defined by Kapranov and Voevodksy~\cite{KV94}, and that 
these should form a monoidal 2-category. Neuchl~\cite{Ne97} studied the 
monoidal 2-categories of representations of Hopf categories in his PhD 
dissertation. In \cite{Ma99} the author of the 
present article studied the categorification of Barrett and Westbury's 
construction and defined spherical 2-categories and the corresponding 
4D state-sum invariants. 

In going from 2D to 3D we had to replace the concept of algebra by that of 
monoidal category, thus allowing for one more layer of structure. Analogously, 
in going from 3D to 4D, we have to add one more layer: besides objects and 
morphisms, we want morphisms between morphisms, which are called 2-morphisms. 
Structures of this sort, called {\it bicategories}, were 
defined by Benabou~\cite{Be67}. He also showed that a 2-category with one object, $X$, can be 
considered as a monoidal category whose objects are the endomorphisms on $X$ 
and whose morphisms are the 2-morphisms between these endomorphisms. The 
tensor product is defined by the composition of the endomorphisms. Note that 
in a bicategory the composition of (1-)morphisms need not to be strictly 
associative: in general there is a non-trivial associator which satisfies the 
pentagon equation. Bicategories are not as exotic as may seem at first. Two 
good examples are the following: the bicategory of all (small) categories and 
the fundamental 2-groupoid of a topological space. The objects of the former 
are all (small) categories, the 1-morphisms are all functors between 
categories, and the 2-morphisms are all natural transformations between the 
functors. In the second example, the objects are the points in the space, the 
1-morphisms are the paths between points, and the 2-morphisms are 
``homotopy classes'' of homotopies between paths. Note that in the first 
example the composition of the 1-morphisms is strictly associative, whereas 
in the second example it is not. Monoidal 2-categories were systematically 
studied by Kapranov and Voevodsky~\cite{KV94}, although some other authors had 
studied particular cases before them. To explain the notion of monoidal 
2-category would take us too far from our main line of reasoning in this 
introduction. Suffice it to mention one important aspect of it, which, 
hopefully, was expected by the reader after reading the earlier paragraphs: 
the 
associator which controls the lack of associativity of the tensor product 
does not satisfy the pentagon equation ``on the nose''. Instead there is a 
{\it modification}, i.e., a natural 2-isomorphism, between the two sides 
of the pentagon equation, called the {\it pentagonator}. This pentagonator 
is required to satisfy a new equation, which is sometimes called the 
{\it non-abelian 4-cocycle relation}. This is completely in conformity with 
the ``basic rule'' of categorification: equations which hold on the nose 
in dimension $n$ are to be substituted by isomorphisms in dimension $n+1$ 
which are required to satisfy new equations. 

\begin{figure}
\centerline{
\epsfbox{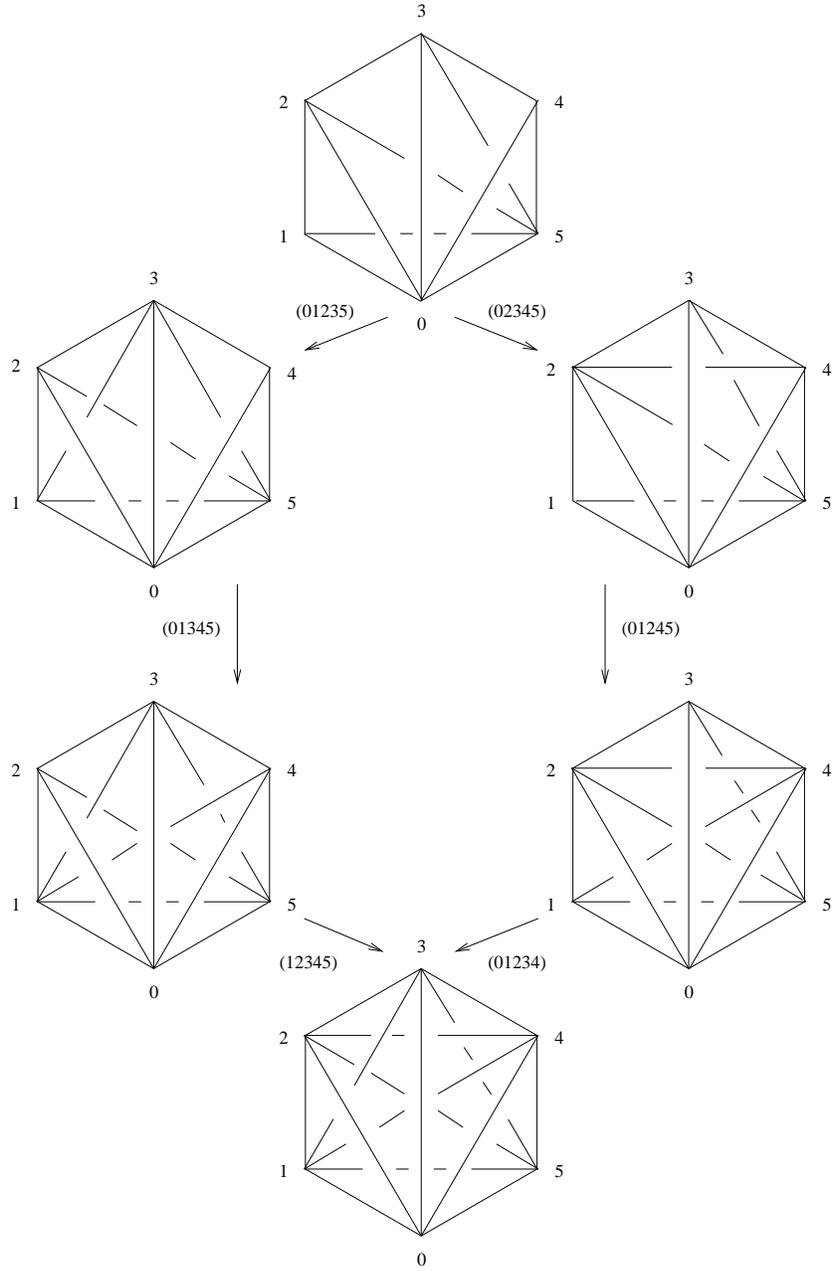}
}
\caption{$3\rightleftharpoons 3$ Pachner move}
\label{pach33}
\end{figure}

\begin{figure}
\centerline{
\epsfbox{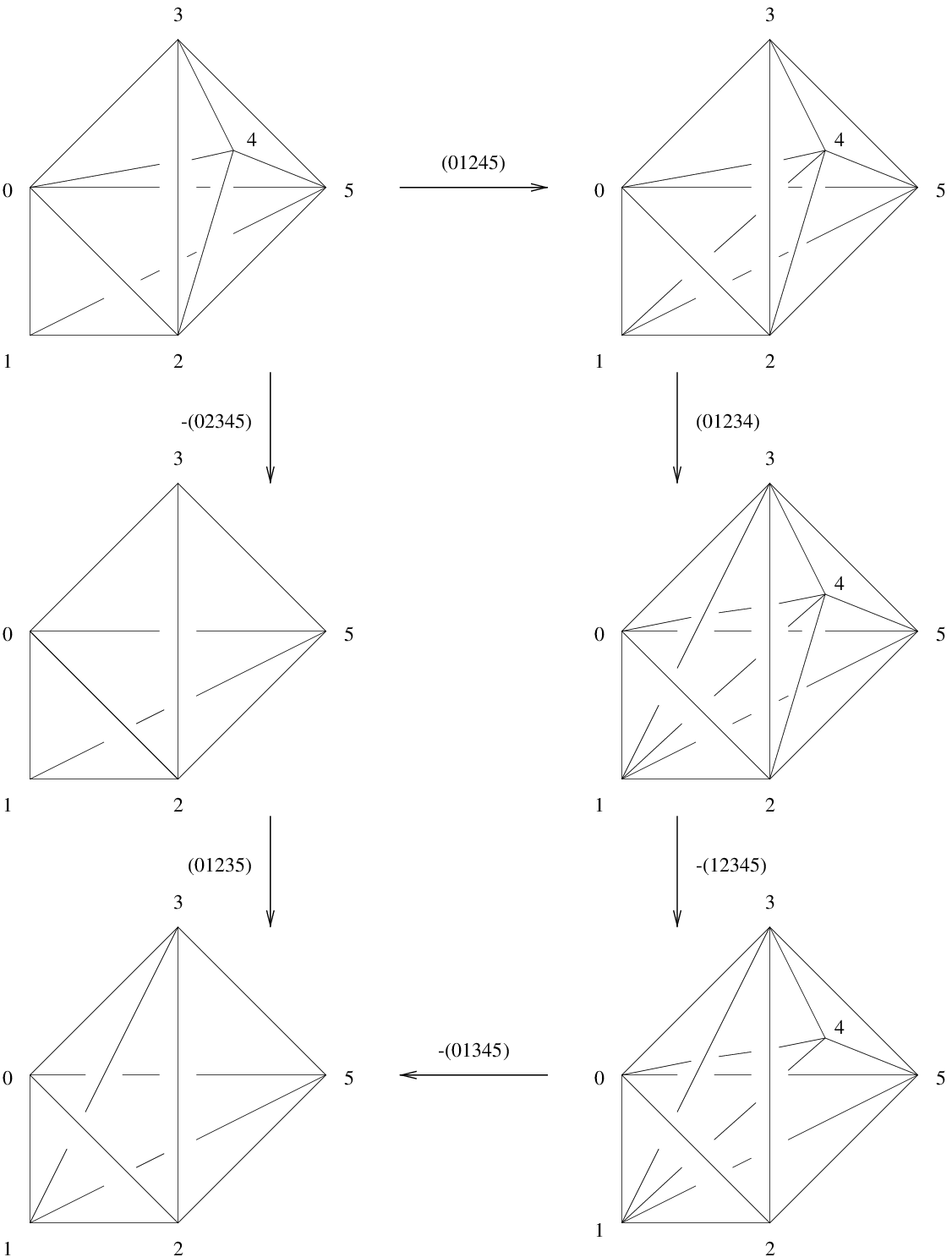}
}
\caption{$2\rightleftharpoons 4$ Pachner move}
\label{pach24}
\end{figure}

\begin{figure}
\centerline{
\epsfbox{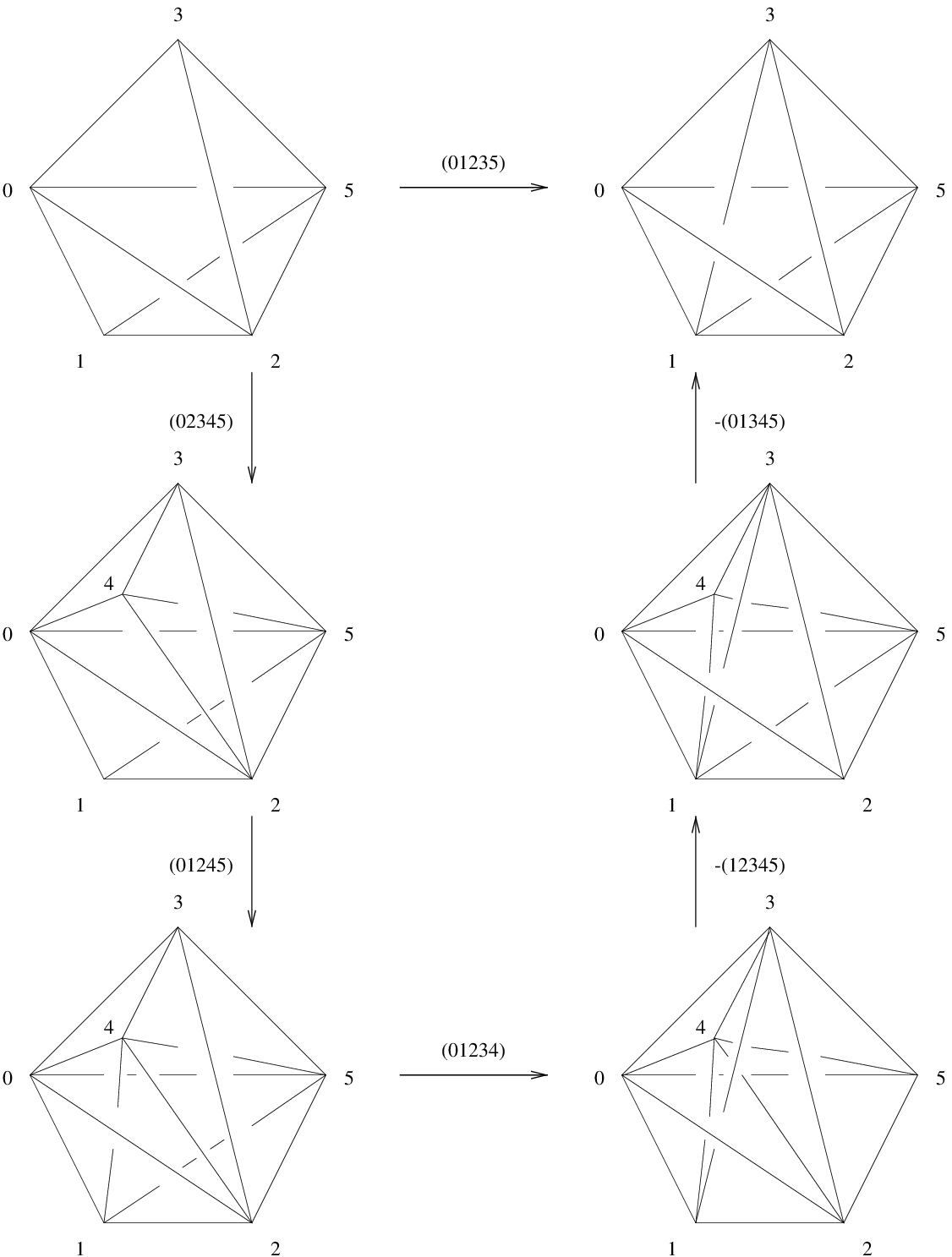}
}
\caption{$1\rightleftharpoons 5$ Pachner move}
\label{pach15}
\end{figure}

For the proof of invariance of the 4D state-sums which we defined in 
\cite{Ma99} it is necessary to 
express the 4D Pachner moves as ``categorifications'' of the 3D Pachner moves, 
which we show in Figs.~\ref{pach33}, \ref{pach24}, \ref{pach15}. Again, these 
figures have been copied from~\cite{CKS97}

Recall that the 3D Pachner moves correspond to partitions of the boundary 
of a 4-simplex into two connected parts with a common boundary. Therefore, 
each arrow in the diagrams of the 4D Pachner moves corresponds to the 
gluing of the two parts of the boundary of a 4-simplex corresponding to 
a 3D Pachner move, after which there is only one way to fill up the missing 
3-cell. In this way one side of a diagram builds up one connected part of 
the boundary of a 5-simplex, whereas the other side builds up the 
complementary part of that boundary.   

In this article we study a particular class of spherical 2-categories and 
the corresponding state-sums. Since in this case the partition function is 
relatively simple, one can prove independence directly, which is 
what we do in Sect.~\ref{Z}. Let us just mention that, as 
in going from 2D to 3D, the algebraic categorification and the 
``topological categorification'' go hand in hand, which ``explains'' the 
invariance of our state-sums. 
 
Conjecturally~\cite{Ma99} the representations of an involutory Hopf category, 
$H$, form a spherical 2-category, $C_H$. Bearing Barrett and Westbury's 
results~\cite{BW95} in mind we conjecture that the Crane-Frenkel invariants 
using $H$ are equal to our invariants using $C_H$. Furthermore, we proved in 
\cite{Ma99} that a spherical 2-category with one object is nothing but 
a tortile category, which is the kind of category that Crane and 
Yetter~\cite{CKY97, CY93} used for their construction of 4-manifold invariants. But for 
a small technical detail, which we explain in Sect.~\ref{N}, it is clear that 
our whole setup generalizes Crane and Yetter's setup. Crane and Yetter's 
invariants are a partial categorification of the Turaev-Viro type 
invariants. Let us finish this part of the introduction by noting that our 
story about categorification is far from complete. We apologize to everyone 
whose contributions to the subject we have not mentioned. In this 
introduction we have tried to give a comprehensive overview of the results 
that lead to our own work, rather than a review of the whole subject. For a 
more complete picture see~\cite{BD982}. 

We summarize our results. For the rest of this paper, let $G$ be any finite group, $H$ any finite abelian group, $R$ any commutative ring with unit and 
involution, which is denoted by $*$, and 
$R^*$ its group of invertible elements. If $r$ is an element in $R$, we 
call $r^*$ its conjugate. In Section~\ref{N} we define 
$\N(G,H,R)$. 
Roughly speaking this is the 2-category of which the objects are finite linear 
combinations of elements of $G$ with non-negative integer coefficients, the 1-endomorphisms of an object $g\in G$ are finite linear combinations of elements of $H$ 
with non-negative integer coefficients, and the 2-endomorphisms of a 1-endormorphism $h\in H$ are elements of the ring $R$. Composition on all levels is 
induced by the group operations, which we write multiplicatively throughout this paper. 

In Section~\ref{M} we define the kind of monoidal structure on $\N(G,H,R)$ 
that we are interested in. We also define an equivalence relation on the 
set of monoidal structures. 

In Section~\ref{Z} we define our state-sum and indicate how we derived its 
definition from our construction in \cite{Ma99}. We do not repeat that 
abstract 
construction here, because it would increase the number of pages considerably and 
might confuse the less category-minded reader. We prove invariance of the state-sums 
that are defined in this paper directly, without going back to the abstract results 
in \cite{Ma99}. For a thorough 
understanding of the results in this paper it is probably better to read \cite{Ma99} 
anyway, but formally the results in this paper are self-contained. 
It is interesting to note that our construction yields the 
``twisted'' version of Yetter's~\cite{Ye93} construction for the case 
of 2-simple path-connected homotopy 2-types. A 2-simple path-connected space 
$E$ is a path-connected space for 
which the action of $\pi_1(E)$, the 
fundamental group, is trivial on $\pi_2(E)$, the second homotopy group. Yetter 
gives a construction of state-sums for arbitrary homotopy 2-types, or 
``categorical groups'', which is equivalent. 

In Section~\ref{P} we relate our results to the theory of homotopy 3-types 
in the form of Postnikov systems. This interpretation in terms of Postnikov 
systems provides a nice link with Freed and 
Quinn's work in \cite{FQ93,Qu98,Qu95}. 

\section{${\mathbb N}(G,H,R)$}
\label{N}

In this section we define the semi-strict monoidal 2-category 
${\mathbb N}(G,H,R)$. We recall that Kapranov and Voevodsky~\cite{KV94} 
defined the more general notion of a weak monoidal 2-category, and that 
Gordon, Power and Street~\cite{GPS95} showed that any weak monoidal 
2-category is equivalent, in an apropriate sense, to a semi-strict monoidal 
2-category. As a matter of fact they proved a more general 
{\it strictification theorem} about 
weak 3-categories, or tricategories, but we only need the case of a 
tricategory with one object which corresponds to a weak monoidal 2-category 
by reindexing of the $k$-morphisms for $0\leq k\leq 3$.  

First let us define the category ${\mathbb N}(H,R)$. 

\begin{definition} ${\mathbb N}(H,R)$ is the $R$-linear finitely semi-simple 
category with the simple objects being precisely 
the elements of $H$, and for which the $R$-module of endomorphisms of an 
object 
$h\in H$ is defined by $\End(h)=R$. The composite of two such 
endomorphisms, $r_1$ and $r_2$, is defined by their product $r_1r_2$ in $R$.
\end{definition}
Note that the objects of ${\mathbb N}(H,R)$ are just finite 
linear combinations 
of elements of $H$ with non-negative integer coefficients. Another way of 
saying this is that the objects are just the elements of the so called 
{\it group rig}~\cite{CY97}, ${\mathbb N}(H)$. If we choose an 
ordering on the elements of $H$, we can represent the morphisms by matrices. 
Let us explain this in a little more detail. Suppose $H$ has order $k$. 
We define the degree of a 
finite linear combination of elements of a 
group with non-negative integer coefficients 
as the sum of the coefficients. We denote the degree of such a 
linear combination $X$ 
by $\hbox{deg}(X)$. A morphism with source $X=n_1h_1+\cdots +n_kh_k$ and 
target $Y=m_1h_1+\cdots +m_kh_k$ can be represented by a 
$\deg(X)\times\deg(Y)$ block diagonal matrix with coefficients in $R$. The $i$-th 
block has size $n_i\times m_i$. Composition is defined by matrix 
multiplication. The product in $H$ induces a monoidal structure on 
${\mathbb N}(H,R)$. Note that we can take ${\mathbb N}(H,R)$ to be symmetric, since 
$H$ is abelian, so $h_1h_2=h_2h_1$. There is also a left duality on 
${\mathbb N}(H,R)$: the left dual of an element $X=n_1h_1+\cdots +n_kh_k$ is 
defined by $X^*=n_1h_1^{-1}+\cdots n_kh_k^{-1}$. The dual of a morphism, 
represented by a matrix, is defined by the conjugate transpose of that matrix. It is not hard to check that this 
symmetry and this duality define a ribbon structure on ${\mathbb N}(H,R)$. 

We are now ready to define the strict 2-category 
${\mathbb N}(G,H,R)$.  
\begin{definition} ${\mathbb N}(G,H,R)$ is the ${\mathbb N}(H,R)$-linear 
finitely semi-simple strict 2-category of which the simple objects are precisely 
the elements of $G$, and for which the ${\mathbb N}(H,R)$-module category of 
endomorphisms on $g\in G$ is defined by $\End(g)={\mathbb N}(H,R)$. 
\end{definition}
Let us explain this definition. The objects of ${\mathbb N}(G,H,R)$ are elements 
of ${\mathbb N}(G)$. Choose an ordering on the 
elements of $G$ and $H$. We can now represent 1- and 2-morphisms by matrices. 
Let $l$ be the order of $G$. A 1-morphism between two objects 
$X=n_1g_1+\ldots + n_lg_l$ and $Y=m_1g_1+\ldots +m_lg_l$ is 
a $\deg(X)\times \deg(Y)$ block diagonal matrix, $f$, with coefficients being 
elements of ${\mathbb N}(H)$. The size of the $i$-th block is equal to 
$n_i\times m_i$. The composition is given by matrix multiplication, where the 
operations on the coefficients are the multiplication and the addition in 
${\mathbb N}(H)$. A 2-morphism between two such 
1-morphisms, $f$ and $g$, is represented by a  
$\deg(X)\times\deg(Y)$ block diagonal matrix, $(\alpha^i_j)$, where  
the coefficient $\alpha^i_j$ is 
a $\deg(f^i_j)\times\deg(g^i_j)$ matrix with coefficients in $R$.  
The horizontal composite 
of two 2-morphisms $\alpha$ and $\beta$, which we denote by $\alpha\circ\beta$, is defined by matrix multiplication, but the operations on the coefficients 
are more complicated than in the case of the 1-morphisms. We define 
$(\alpha\circ\beta)^i_j=\oplus_{k}(\alpha^i_k\otimes\beta^k_j)$. Note that the coefficients $\alpha^i_k$ and $\beta^k_j$ are matrices 
themselves with coefficients in $R$. In general the tensor product of two matrices, 
$x$ and $y$, is defined by $(x\otimes y)^{ij}_{kl}=x^i_ky^j_l$, and the direct sum is 
defined by 
$$x\oplus y=\begin{pmatrix}x&0\\ 0&y\end{pmatrix}.$$ The vertical composite of two 2-morphisms, 
$\alpha$ and $\beta$, which we denote by $\alpha\cdot\beta$, is defined by 
coefficientwise multiplication, i.e., 
$(\alpha\cdot\beta)^i_j=\alpha^i_j\beta^i_j$. Note that, just as in 
the completely coordinatized version of the monoidal 
2-category of 2-vector spaces, $2\Vect_{cc}$, defined by Kapranov and 
Voevodsky~\cite{KV94}, the coefficients of the 2-morphisms are matrices 
themselves, so their 
multiplication 
is given by matrix multiplication. Note also that this multiplication is well 
defined for any pair of composable 2-morphisms: for any $\alpha\colon f\to g$ and 
$\beta\colon g\to h$, the matrix $\alpha^i_j$ has size $f^i_j\times g^i_j$ and 
$\beta^i_j$ has size $g^i_j\times h^i_j$. As always, we write the composites 
of 1- and 2-morphisms in the 
diagrammatic order. It is easy to check that all compositions are 
strictly associative. 

The semi-strict monoidal structure on ${\mathbb N}(G,H,R)$ 
is induced by the multiplication in $G$, $H$, and $R$. For any two objects 
$X$ and $Y$, we define $X\otimes Y=XY$. For any 1-morphism 
$f$ and any object $Y$ we define $f\otimes Y=f\otimes 1_Y$, where $1_Y$ is 
the identity on $Y$. In terms of 
coefficients this becomes $(f\otimes 1_Y)^{ij}_{kl}=f^i_k\delta^j_l$, with $\delta$ 
being the Kronecker delta. For any 2-morphism $\alpha$ and any object $Y$, 
we define $\alpha\otimes Y=\alpha\otimes 1_{1_Y}$, where $1_{1_Y}$ is the 
identity 2-morphism of the identity 1-morphism on $Y$. In terms of 
coefficients this becomes $(\alpha\times 1_{1_Y})^{ij}_{kl}=
\alpha^i_k\otimes (1_Y)^j_l$. Note that $\alpha$ and $1_Y$ are matrices themselves, so that 
there tensor product differs from their product in general. Likewise we define 
$Y\otimes f$ and $Y\otimes\alpha$. It is easy to check that these tensor 
products are strictly associative. The tensorator 
$$\bigotimes\nolimits_{f,g}\colon 
(f\otimes Y)(X'\otimes g)\Rightarrow (X\otimes g)(f\otimes Y')$$
of two 1-morphisms $f\colon X\to X'$ and $g\colon Y\to Y'$ is defined by the 
operator that interchanges the two tensor factors. Concretely, we have 
$$[(f\otimes Y)(X'\otimes g)]^{ij}_{kl}=f^i_kg^j_l$$
and
$$[(X\otimes g)(f\otimes Y')]^{ij}_{kl}=g^i_kf^j_l,$$ 
so the tensorator becomes 
$$(\bigotimes\nolimits_{f,g})^{ij}_{kl}=P_{f^i_kg^j_l},$$ 
where $P_{f^i_kg^j_l}$ is the $\deg(f^i_k)\times\deg(g^j_l)$ matrix with coefficients 
$(P_{f^i_kg^j_l})^{mn}_{rs}=\delta^m_s\delta^n_r$. This defines a semi-strict 
monoidal structure on ${\mathbb N}(G,H,R)$. We remark that 
${\mathbb N}(G,\{1\},{\mathbb C})$ is equal to $2\Hilb[G]$, which we defined in 
\cite{Ma99}. If $G=\{1\}$ 
also, then we recover the definition of the completely coordinatized version of the monoidal 
2-category of 2-vector spaces, $2\Vect_{cc}$, see~\cite{KV94}. 

We can now define the left-duality on ${\mathbb N}(G,H,R)$ in terms of matrices. 
The dual of an object $X=n_1g_1+\cdots +n_lg_l$ is defined by 
$X^*=n_1g^{-1}+\cdots +n_lg^{-1}$. The dual of a 1-morphism, represented by 
a matrix $f$, is defined by the transpose of $f$ with dual coefficients. 
Recall that the dual of a coefficient $m_1y_1+\cdots+m_kh_k$ is defined by 
$m_1h_1^{-1}+\cdots m_kh_k^{-1}$. The 
dual of a 2-morphism, $\alpha$, is defined by the matrix $\alpha^*$, where 
$(\alpha^*)^i_j$ is the conjugate transpose of $\alpha^i_j$. This definition 
of the duality the author derived from the definition of duality in 
the monoidal 2-category of 2-Hilbert spaces by Baez~\cite{Ba97}. Baez defined 
a ``weak'' version of 2-Hilbert spaces, the underlying monoidal 2-category 
of which is equivalent to $2\Vect$, the non-coordinatized version of 2-vector 
spaces~\cite{KV94}. For our application it is better to work with the 
completely coordinatized version, $2\Hilb_{cc}$, which is 
semi-strict~\cite{Ma99}. Just as for $2\Hilb_{cc}$, it is not 
hard to show that this duality satisfies the {\it spherical 
conditions}~\cite{Ma99}. We do not prove this, because we do not 
need it explicitly in this article.   

\section{Monoidal structures on ${\mathbb N}(G,H,R)$}
\label{M}

The most interesting state-sums are related to weakenings of the 
semi-strict monoidal structure on ${\mathbb N}(G,H,R)$. 
These weakenings can be defined, essentially, by following the definition of 
a monoidal 2-category by Kapranov and Voevodsky \cite{KV94}. We repeat that 
Kapranov and Voevodsky's definition coincides with that of Gordon, Power and 
Street's~\cite{GPS95} definition of a tricategory with one object. Kapranov 
and Voevodsky, using MacLane and Pare's coherence theorem~\cite{MP85}, assume 
that the underlying 2-category is 
strict. We do not want to make this assumption, because it is too 
restrictive for our purpose. Therefore 
we have to keep in account the non-associativity of the composition of 
the edges in the diagrams in \cite{KV94}. We assume that 
this non-associativity is controlled by a {\it coherent associator}, so it 
does not 
matter how we 
choose to parenthesize the boundary 1-morphisms in the diagrams. We just make 
one choice 
and work out the diagrams. Any other choice will lead to equivalent diagrams. 
Before going on, let us have a look at this {\it associator}. Note that 
we can restrict our attention to the 1-morphisms in $\End(1)$, because the 
general case then follows by linearity. An associator on 1-morphisms in 
$\End(1)$ is a family 
of 2-isomorphisms 
$$\alpha^1_{h_1,h_2,h_3}\colon h_1(h_2h_3)\To (h_1h_2)h_3,$$ 
indexed by triples of 1-morphisms. Since all 1-morphisms are sums of simple 
1-morphisms, which are simply elements of $H$ in this case, we only have 
to define $\alpha^1$ on triples of simple 1-morphisms; the general definition 
then follows by extending $\alpha^1$ linearly. Note that $h_1(h_2h_3)$ and 
$(h_1h_2)h_3$ are the same 1-morphism. The associator is a natural isomorphism 
between the two functors $\End(1)\times\End(1)\times\End(1)\to\End(1)$ 
which define the two different ways of composing three 1-morphisms, indicated 
by the different bracketings. Hence $\alpha^1_{h_1,h_2,h_3}$ is just a 
2-automorphism on $h_1h_2h_3$, i.e., an element of $R^*$. Thus 
we can define $\alpha^1$ as a function $H\times H\times H\to R^*$. As 
usual in category theory, we have to impose a condition on $\alpha^1$, 
called a {\it coherence relation}, in order to maintain control over the 
bracketing. There are five different ways of composing four 1-morphisms, and 
two different ways of rebracketing the composite going from 
right-to-left bracketing to left-to-right bracketing. This leads to the 
{\it pentagon diagram} which we showed in the introduction and which we 
repeat here.
$$
\begin{CD}
f(g(hj))@>\alpha^1_{f,g,hj}>> (fg)(hj)@>\alpha^1_{fg,h,j}>> ((fg)h)j\\
@V\alpha^1_{g,h,j}VV      &           &       @AA\alpha^1_{f,g,h}A\\
f((gh)j)  &    @>\alpha^1_{f,gh,j}>>   &      (f(gh))j
\end{CD}
$$
MacLane and Pare's coherence theorem~\cite{MP85} implies that, if 
we require the pentagon diagram to be commutative, that any two strings of 
$\alpha^1$'s, i.e., composites of an arbitrary number of associators, with 
the same source and target are equal. This means, for example, that two 
different 
algorithms that rebracket composites of 1-morphisms will always end up 
using the same 2-isomorphism, although one algorithm may use a different 
decomposition of this 2-isomorphism than the other. This is why we impose 
the condition
$$\alpha^1_{h_2,h_3,h_4}
(\alpha^1_{h_1h_2,h_3,h_4})^{-1}\alpha^1_{h_1,h_2h_3,h_4}
(\alpha^1_{h_1,h_2,h_3h_4})^{-1}\alpha^1_{h_1,h_2,h_3}=1,$$
on $\alpha^1$. The maps and conditions in Def.~\ref{sweakmonstruc} are 
obtained from Kapranov and Voevodsky's diagrams in a similar way. 
In the 
following definition we use their hieroglyphic notation to 
indicate these 
diagrams. 

\begin{definition}
\label{sweakmonstruc}
A {\rm semi-weak monoidal 2-category structure} on ${\mathbb N}(G,H,R)$ consists 
of the following maps:
\begin{description}
\item[0-associator] $\alpha^0\colon G\times G\times G\to H$, which corresponds 
to a family of simple invertible 1-morphisms $\alpha^0_{g_1,g_2,g_3}\colon 
g_1(g_2g_3)\to (g_1g_2)g_3$.
\item[pentagonator] $\pi\colon G\times G\times G\times G\to R^*$, which 
corresponds to a family of invertible 2-morphisms 
\begin{equation*}
\begin{CD}
g_1(g_2(g_3g_4))@>\alpha^0_{g_1,g_2,g_3g_4}>> (g_1g_2)(g_3g_4) @>\alpha^0_{g_1g_2,g_3,g_4}>> ((g_1g_2)g_3)g_4\\
@V\alpha^0_{g_2,g_3,g_4}VV \Big\Downarrow\pi_{g_1,g_2,g_3,g_4}&&
@AA\alpha^0_{g_1,g_2,g_3}A\\
g_1((g_2g_3)g_4)&@>\alpha^0_{g_1,g_2g_3,g_4}>> &(g_1(g_2g_3))g_4
\end{CD}
\end{equation*}
\item[1-associator] $\alpha^1\colon H\times H\times H\to R^*$, which 
corresponds to a family of invertible 2-morphisms $\alpha^1_{h_1,h_2,h_3}
\colon h_1(h_2h_3)\To (h_1h_2)h_3$. 
\item[tensorator] $\tau\colon H\times H\to R^*$, which corresponds to 
a family of invertible 2-morphisms 
\begin{equation*}
\begin{CD}
g_1g_2 & @>h_1>> & g_1g_2\\
@Vh_2VV \Big\Downarrow\tau_{h_1,h_2} &&@VVh_2V\\
g_1g_2 & @>h_1>> & g_1g_2
\end{CD}
\end{equation*}
\item[interchanger1] $\iota^1\colon H\times G\times G\to R^*$, which 
corresponds to a family of invertible 2-morphisms
\begin{equation*}
\begin{CD}
g_1(g_2g_3) & @>\alpha^0_{g_1,g_2,g_3}>> & (g_1g_2)g_3\\
@Vh_1VV \Big\Downarrow\iota^1_{h_1,g_2,g_3} && @VVh_1V\\
g_1(g_2g_3) & @>\alpha^0_{g_1,g_2,g_3}>> & (g_1g_2)g_3
\end{CD}
\end{equation*} 
\item[interchanger2] $\iota^2\colon G\times H\times G\to R^*$, which 
corresponds to a family of invertible 2-morphisms
\begin{equation*}
\begin{CD}
g_1(g_2g_3) & @>\alpha^0_{g_1,g_2,g_3}>> & (g_1g_2)g_3\\
@Vh_2VV \Big\Downarrow\iota^2_{g_1,h_2,g_3} && @VVh_2V\\
g_1(g_2g_3) & @>\alpha^0_{g_1,g_2,g_3}>> & (g_1g_2)g_3
\end{CD}
\end{equation*} 
\item[interchanger3] $\iota^3\colon G\times G\times H\to R^*$, which 
corresponds to a family of invertible 2-morphisms
\begin{equation*}
\begin{CD}
g_1(g_2g_3) & @>\alpha^0_{g_1,g_2,g_3}>> & (g_1g_2)g_3\\
@Vh_3VV \Big\Downarrow\iota^3_{g_1,g_2,h_3} && @VVh_3V\\
g_1(g_2g_3) & @>\alpha^0_{g_1,g_2,g_3}>> & (g_1g_2)g_3
\end{CD}
\end{equation*} 
\end{description}
All these maps are required to be normalized, i.e., their value is equal to 
$1$ whenever 
one of the factors of their argument is equal to $1$. 
Furthermore these maps are required to satisfy the following identities:
\begin{description}\label{weakmon}
\item[$(\bullet\otimes\bullet\otimes\bullet\otimes\bullet)$] 
$$\alpha^0_{g_2,g_3,g_4}(\alpha^0_{g_1g_2,g_3,g_4})^{-1}\alpha^0_{g_1,g_2g_3,g_4}
(\alpha^0_{g_1,g_2,g_3g_4})^{-1}\alpha^0_{g_1,g_2,g_3}=1.$$
\item[$(\rightarrow\rightarrow\rightarrow\rightarrow)$] 
$$\alpha^1_{h_2,h_3,h_4}
(\alpha^1_{h_1h_2,h_3,h_4})^{-1}\alpha^1_{h_1,h_2h_3,h_4}
(\alpha^1_{h_1,h_2,h_3h_4})^{-1}\alpha^1_{h_1,h_2,h_3}=1.$$
\end{description}
In the following identities we avoid writing $\alpha^1$ constantly and bracket the 
remaining maps with $\lceil\rceil$ following Crane and Yetter's~\cite{CY972} notation. As 
explained in \cite{CY972} this notation means that the source and target 1-morphisms are assumed to be parenthesized from left to right. The brackets denote 
the strings of 1-associators that are required to make the 2-morphisms 
composable under this assumption. The usage of these brackets is unambiguous 
by the coherence relation of the 
1-associator, which corresponds to $(\rightarrow\rightarrow\rightarrow\rightarrow)$. 
\begin{description}
\item[$(\bullet\otimes\bullet\otimes\bullet)$]
$$\lceil\tau_{h_1h_2,h_3}\rceil=\lceil\tau_{h_2,h_3}\rceil\lceil\tau_{h_1,h_3}\rceil,$$
and
$$\lceil\tau_{h_1,h_2h_3}\rceil=\lceil\tau_{h_1,h_2}\rceil\lceil\tau_{h_1,h_3}\rceil.$$
\item[$(\bullet\otimes\bullet\otimes\bullet\otimes\bullet\otimes\bullet)$]
$$\lceil\pi_{g_2,g_3,g_4,g_5}\rceil\lceil\pi_{g_1,g_2g_3,g_4,g_5}\rceil\lceil
\iota^1_{\alpha^0_{g_1,g_2,g_3};g_4;g_5}\rceil\lceil\pi_{g_1,g_2,g_3,g_4g_5}\rceil=$$
$$\lceil\iota^2_{g_1;\alpha^0_{g_2,g_3,g_4};g_5}\rceil\lceil\pi_{g_1,g_2,g_3,g_4}\rceil 
\lceil\pi_{g_1,g_2,g_3g_4,g_5}\rceil\lceil(\iota^3_{g_1;g_2;\alpha^0_{g_3,g_4,g_5}})^{-1}\rceil\lceil\pi_{g_1g_2,g_3,g_4,g_5}\rceil.$$
\item[$(\rightarrow\otimes\bullet\otimes\bullet\otimes\bullet)$]
$$\lceil\iota^1_{h,g_2,g_3}\rceil\lceil\iota^1_{h,g_2g_3,g_4}\rceil\lceil\tau^{-1}_{h;\alpha^0
_{g_2,g_3,g_4}}\rceil=\lceil\iota^1_{h,g_3,g_4}\rceil\lceil\iota^1_{h,g_2,g_3g_4}\rceil.$$
\item[$(\bullet\otimes\rightarrow\otimes\bullet\otimes\bullet)$]
$$\lceil\iota^2_{g_1,h,g_3g_4}\rceil=\lceil\iota^2_{g_1,h,g_3}\rceil
\lceil\iota^2_{g_1,h,g_4}\rceil.$$
\item[$(\bullet\otimes\bullet\otimes\rightarrow\otimes\bullet)$]
$$\lceil\iota^2_{g_1g_2,h,g_4}\rceil=\lceil\iota^2_{g_1,h,g_4}\rceil
\lceil\iota^2_{g_2,h,g_4}\rceil.$$
\item[$(\bullet\otimes\bullet\otimes\bullet\otimes\rightarrow)$]
$$\lceil\tau_{\alpha^0_{g_1,g_2,g_3};h}\rceil\lceil\iota^3_{g_1,g_2g_3,h}
\rceil\lceil\iota^3_{g_2,g_3,h}\rceil=\lceil\iota^3_{g_1g_2,g_3,h}\rceil
\lceil\iota^3_{g_1,g_2,h}\rceil.$$
\item[$(\rightarrow\rightarrow\otimes\bullet\otimes\bullet)$]
$$\lceil\iota^1_{h_1h_2,
g_2,g_3}\rceil=\lceil\iota^1_{h_1,g_2,g_3}\rceil\lceil\iota^1_{h_2,g_2,g_3}
\rceil.$$
\item[$(\bullet\otimes\rightarrow\rightarrow\otimes\bullet)$]
$$\lceil\iota^2_{g_1,h_1h_2,g_3}\rceil=\lceil\iota^2_{g_1,h_1,g_3}\rceil
\lceil\iota^2_{g_1,h_2,g_3}\rceil.$$
\item[$(\bullet\otimes\bullet\otimes\rightarrow\rightarrow)$]
$$\lceil\iota^3_{g_1,g_2,h_1h_2}\rceil=\lceil\iota^3_{g_1,g_2,h_1}\rceil
\lceil\iota^3_{g_1,g_2,h_2}\rceil.$$
\end{description}
\end{definition}
Let us briefly comment on these maps and relations. There is only one 
structural 1-morphism: the 0-associator, $\alpha^0$. It controls 
the non-associativity of the tensor product on objects and is given by a 
family of invertible 1-morphisms indexed by triples of objects. 
It suffices to define $\alpha^0$ on 
simple objects, i.e., elements of $G$. We assume in our definition that all 
structural 
1-morphisms are simple. Therefore we define $\alpha^0$ to take 
values in $H$. It is now easy to 
derive the 3-cocycle condition in 
$(\bullet\otimes\bullet\otimes\bullet\otimes\bullet)$ from the corresponding 
diagram in \cite{KV94}. Note that this diagram is just the pentagon diagram 
for objects and 1-morphisms. All other maps in Definition~\ref{sweakmonstruc} 
are structural 2-morphisms. Since they are also assumed to be invertible, 
they take values in $R^*$. It suffices 
to index them by simple objects, i.e., elements in $G$, and simple 
1-morphisms, i.e., 
elements in $H$. The list of maps and relations now follows easily from 
Kapranov and Voevodsky's definitions. The pentagonator, $\pi$, controls the 
non-commutativity of the pentagon diagram for the 0-associator. This pentagon 
diagram corresponds to $(\bullet\otimes\bullet\otimes\bullet\otimes\bullet)$. 
As we already explained, the 1-associator, $\alpha^1$, 
controls the non-associativity of the composition of the 1-morphisms. The 
tensorator, $\tau$, is a weakening of the tensorator in the semi-strict 
monoidal structure on ${\mathbb N}(G,H,R)$. Finally, the interchangers, 
$\iota^i$ for $i=1,2,3$, define the pseudo-naturality of $\alpha^0$. 
All relations are coherence relations which ensure that the composites of any 
two strings 
of structural maps with the same source and target are equal. The 
assumption that all morphisms are simple is 
restrictive, but is inspired by the relation with homotopy theory, as 
explained in Section~\ref{P}. 
A second reason for this assumption is that the calculations, which are 
not easy anyway, become much simpler under this assumption.   
We call these structures semi-weak, because we assume the units to be strict 
and the tensor product of an object with a 1- or 2-morphism to be trivial. 
Therefore some of the structural 1- and 2-morphisms in \cite{KV94} become 
identities. This also explains why we have fewer coherence relations than 
Kapranov and Voevodsky have in \cite{KV94}. Note that $(\bullet\otimes\bullet
\otimes\bullet\otimes\bullet)$ and $(\rightarrow\rightarrow\rightarrow
\rightarrow)$ are 3-cocycle conditions. The relations in $(\bullet\otimes\bullet
\otimes\bullet)$ are called the {\it hexagon relations} and together with the 
3-cocycle relation in $(\rightarrow\rightarrow\rightarrow\rightarrow)$ they 
define the structure of a braided monoidal category on $\End(1)$, 
see~\cite{Ma99}. 
The coherence cube $(\rightarrow
\otimes\bullet\otimes\bullet)$ in Kapranov and Voevodsky's paper becomes a 
consequence of the hexagon relations and the triviality of the tensor 
product of a simple object with a 1- or 2-morphism in our setup.  

For this particular class of monoidal 2-categories it is easy to define when 
they are ``equivalent''. We follow Gordon, Power and Street's~\cite{GPS95} definition of 
{\it triequivalence} of tricategories for the special case of 
tricategories with one object, which can be considered as weak monoidal 
2-categories. 

\begin{definition}
\label{equival}
We say that two semi-weak monoidal 2-category structures on 
${\mathbb N}(G,H,R)$, as defined in Def.~\ref{sweakmonstruc}, are 
{\rm 2-equivalent} if there exist
\begin{enumerate}
\item Automorphisms $G\to G$, $H\to H$, and $R\to R$, which we denote 
by $g\mapsto \ol{g}$, $h\mapsto \ol{h}$, and $r\mapsto \ol{r}$. The first 
two automorphisms are required to be group automorphisms, the third one is 
required to be a unital ring automorphism which preserves the involution. 
\item A map $\mu\colon H\times H\to R^*$, which corresponds to a family of 
invertible 2-morphisms 
$\mu_{h_1,h_2}\colon \ol{h_1h_2}\To\ol{h}_1\ol{h}_2$.
\item A map $\Phi\colon G\times G\to H$, which corresponds to a family of 
simple invertible 1-morphisms $\Phi_{g_1,g_2}\colon \ol{g_1g_2}\to 
\ol{g}_1\ol{g}_2$. 
\item A map $\phi\colon G\times G\times G\to R^*$, which corresponds to 
a family of invertible 2-morphisms
\begin{equation*}
\begin{CD}
\ol{g_1(g_2g_3)}@>\Phi_{g_1,g_2g_3}>>\ol{g}_1(\ol{g_2g_3})@>\Phi_{g_2,g_3}>>
\ol{g}_1(\ol{g}_2\ol{g}_3)\\
@V\ol{\alpha^0_{g_1,g_2,g_3}}VV \Big\Uparrow\phi_{g_1,g_2,g_3} &&@VV(\alpha^0)'_{\ol{g}_1,\ol{g}_2,\ol{g}_3}V\\
\ol{(g_1g_2)g_3}@>\Phi_{g_1g_2,g_3}>> (\ol{g_1g_2})\ol{g}_3 @>\Phi_{g_1,g_2}>>
(\ol{g}_1\ol{g}_2)\ol{g}_3
\end{CD}
\end{equation*}
\item A map $\psi\colon H\times G\to R^*$, which corresponds to a family of 
invertible 2-morphisms
\begin{equation*}
\begin{CD}
\ol{g_1g_2}&@>\Phi_{g_1,g_2}>>& \ol{g}_1\ol{g}_2\\
@V\ol{h}_1VV \Big\Downarrow\psi_{h_1,g_2}&& @VV\ol{h}_1V\\
\ol{g_1g_2}&@>\Phi_{g_1,g_2}>>&\ol{g}_1\ol{g}_2
\end{CD}
\end{equation*}
\item A map $\chi\colon G\times H\to R^*$, which corresponds to a family of 
invertible 2-morphisms
\begin{equation*}
\begin{CD}
\ol{g_1g_2}&@>\Phi_{g_1,g_2}>>& \ol{g}_1\ol{g}_2\\
@V\ol{h}_2VV \Big\Downarrow\chi_{g_1,h_2}&& @VV\ol{h}_2V\\
\ol{g_1g_2}&@>\Phi_{g_1,g_2}>>&\ol{g}_1\ol{g}_2
\end{CD}
\end{equation*}
\end{enumerate} 
All these maps are required to be normalized.
Furthermore, they should satisfy
\begin{enumerate}
\item $$\ol{\alpha^0_{g_1,g_2,g_3}}\Phi_{g_1g_2,g_3}\Phi_{g_1,g_2}=
\Phi_{g_1,g_2g_3}\Phi_{g_2,g_3}(\alpha^0)'_{\ol{g}_1,\ol{g}_2,\ol{g}_3}.$$
\item $$\ol{\alpha^1_{h_1,h_2,h_3}}\mu_{h_1h_2,h_3}\mu_{h_1,h_2}=
\mu_{h_1,h_2h_3}\mu_{h_2,h_3}(\alpha^1)'_{\ol{h}_1,\ol{h}_2,\ol{h}_3}.$$
From now on we do not write $\alpha^1$, $(\alpha^1)'$, or 
$\mu$ any longer. 
As explained in Def.~\ref{sweakmonstruc}, we use the brackets $\lceil\rceil$, 
which ``absorb'' these three maps.
\item $$\lceil\overline{\tau_{h_1,h_2}}\rceil=\lceil 
\tau'_{\ol{h}_1,\ol{h}_2}\rceil.$$ 
\item $$\lceil\psi_{h,g_2}\rceil\lceil\psi_{h,g_3}\rceil\lceil
\ol{\iota^1_{h,g_2,g_3}}\rceil=\lceil(\iota^1)'_{\ol{h},\ol{g}_2,\ol{g}_3}
\lceil(\tau')^{-1}_{\ol{h},\Phi_{g_2,g_3}}\rceil\lceil
\psi_{h,g_2g_3}\rceil.$$
\item $$\lceil\ol{\iota^2_{g_1,h,g_3}}\rceil=\lceil(\iota^2)'_{\ol{g}_1,
\ol{h},\ol{g}_3}\rceil.$$
\item $$\lceil \tau'_{\Phi_{g_1,g_2},\ol{h}}\rceil\lceil\chi_{g_1g_2,h}
\rceil\lceil\ol{\iota^3_{g_1,g_2,h}}\rceil=\lceil(\iota^3)'_{\ol{g}_1,
\ol{g}_2,\ol{h}}\rceil\lceil\chi_{g_2,h}\rceil\lceil\chi_{g_1,h}\rceil.$$
\item $$\lceil\psi_{h_1h_2,g}\rceil=\lceil\psi_{h_1,g}\rceil\lceil\psi_{h_2,g
}\rceil.$$
\item $$\lceil\chi_{g,h_1h_2}\rceil=\lceil\chi_{g,h_1}\rceil\lceil\chi_{g,h_2}
\rceil.$$
\item 
\begin{equation*}
\begin{split}
&\lceil\psi^{-1}_{\ol{\alpha^0_{g_1,g_2,g_3}},g_4}\rceil\lceil
\phi_{g_1,g_2,g_3}\rceil\lceil\phi_{g_1,g_2g_3,g_4}\rceil\lceil
(\iota^2)'_{\ol{g}_1,\Phi_{g_2,g_3},\ol{g}_4}\rceil \\
&\times\lceil
\chi^{-1}_{g_1,\alpha^0_{g_2,g_3,g_4}}\rceil\lceil\phi_{g_2,g_3,g_4}\rceil
\lceil\pi'_{\ol{g}_1,\ol{g}_2,\ol{g}_3,\ol{g}_4}\rceil \\
=\ &
\lceil\ol{\pi_{g_1,g_2,g_3,g_4}}\rceil\lceil\phi_{g_1g_2,g_3,g_4}\rceil\lceil
(\iota^1)'_{\Phi_{g_1,g_2},\ol{g}_3,\ol{g}_4}\rceil\lceil
(\tau')^{-1}_{\Phi_{g_1,g_2},\Phi_{g_3,g_4}}\rceil \\
&\times\lceil
\phi_{g_1,g_2,g_3g_4}\rceil\lceil
(\iota^3)'_{\ol{g}_1,\ol{g}_2,\Phi_{g_3,g_4}}\rceil.
\end{split}
\end{equation*}
\end{enumerate}
\end{definition}   
Again, writing down the diagrams makes the conditions in Def.~\ref{equival} 
more comprehensible. The diagrams corresponding to the first seven conditions 
follow easily from the formulas, the diagram corresponding to condition 
8 can be found in \cite{GPS95}. Def.~\ref{equival} defines an equivalence 
relation on the semi-weak monoidal 2-category structures on 
${\mathbb N}(G,H,R)$. 

The duality on the semi-strict ${\mathbb N}(G,H,R)$ is compatible with any 
semi-weak monoidal structure. Note that, by definition, all structural 
2-morphisms are taken to be {\it unital}. Recall that a 2-morphism, $\alpha$, is 
called unital if it is 
invertible and if its dual equals its inverse. 

At the end of the next section we give some examples of semi-weak monoidal 
structures on ${\mathbb N}(G,H,R)$, for $G=\{1\}$, $H={\mathbb Z}/p{\mathbb Z}$, 
$R={\mathbb C}$, and $G=H={\mathbb Z}/p{\mathbb Z}$, $R={\mathbb C}$, 
respectively. These examples are due to Birmingham and Rakowski~\cite{BR94,
BR95,BR96}. Since they also did some calculations of the related state-sums, 
we prefer to explain their results, which fit nicely into our setup, after 
defining our more general state-sums and showing that they are invariant.
 
\section{The state-sums}
\label{Z}
\begin{figure}
\centerline{
\epsfbox{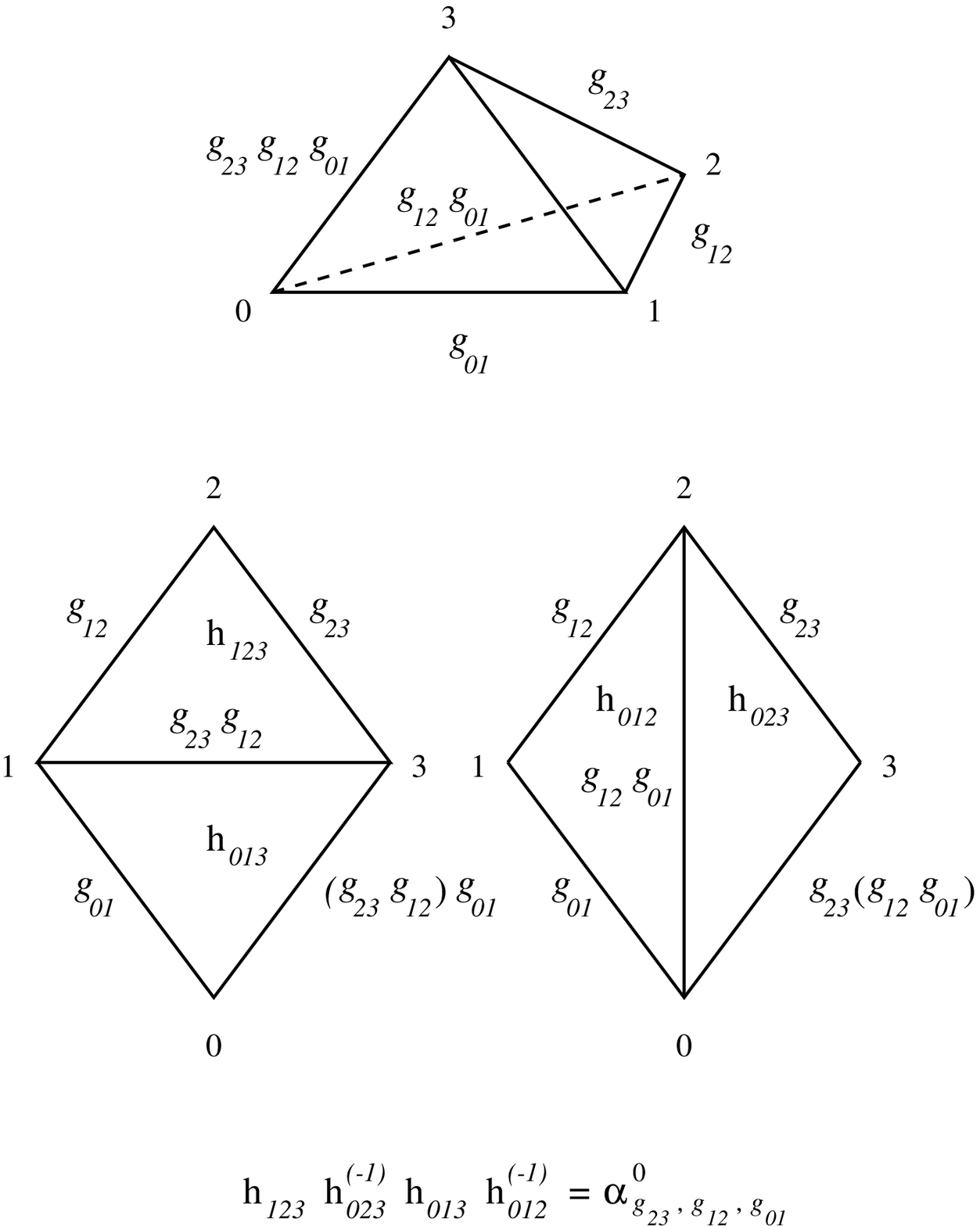}
}
\caption{labellings}
\label{label}
\end{figure}

Fix a semi-weak monoidal 
structure on ${\mathbb N}(G,H,R)$. Henceforth a 4-manifold means a closed 
oriented PL manifold of dimension 4 and any triangulation is assumed to have 
a total ordering on its vertices. Let $M$ be a 4-manifold and ${\cal T}$ a 
triangulation of $M$. Following our setup in \cite{Ma99}, we 
label the edges of ${\cal T}$ with 
elements of $G$ and label the faces, i.e., triangles, with elements of 
$H$. If $(ijk)$ is a face in ${\cal T}$, then we impose the 
condition 
$$g_{ij}g_{ik}^{-1}g_{jk}=1\in G$$
on the labels of the edges. If $(ijkl)$ is 
a 3-simplex in ${\cal T}$, then we require the 
condition
$$h_{jkl}h_{ikl}^{-1}h_{ijl}h_{ijk}^{-1}=\alpha^0_{g_{kl},g_{jk},g_{ij}}\in H$$
to hold true. 
We call these conditions the 
'local semi-flatness' 
conditions and Fig.~\ref{label} shows them diagrammatically. 

These conditions follow naturally from the general setup in~\cite{Ma99}, 
because labellings that do not satisfy these conditions 
correspond to zero terms in the state sum. For example, if 
$$g_{ij}g_{ik}^{-1}g_{jk}\neq 1,$$
then $$\Hom(g_{ik},g_{jk}g_{ij})=\{0\}.$$
We now define the partition function 
on any 4-simplex, $(ijklm)$.
\begin{definition}
\label{Z((ijklm))}
\begin{equation*}
\begin{split}
Z((ijklm))=\ &\lceil (\iota^3_{g_{lm},g_{kl},h_{ijk}})^{-1}\rceil\lceil
\tau_{h_{klm},h_{ijk}}\rceil\lceil(\iota^1_
{h_{klm},g_{jk},g_{ij}})^{-1}\rceil\\
&\times\lceil
\iota^2_{g_{lm},h_{jkl},g_{ij}}\rceil
\lceil\pi_{g_{lm},g_{kl},g_{jk},g_{ij}}\rceil.
\end{split}
\end{equation*}
\end{definition}
We derived this partition function from the abstract one defined 
in~\cite{Ma99} by keeping track of the parentheses around the 
objects and the 1-morphisms that are involved. Note that the 
'funny brackets' $\lceil\rceil$ are very 
helpful here; without them the definition of the partition function would 
contain at least 22 factors. Note also that $h_3$ does not show up in our 
notation 
of the partition function; it is hidden by the brackets. The advantages of 
our notation for the proof of invariance of our state-sum outweighs this 
minor drawback. 

There are two special cases: 
\begin{enumerate}
\item \#H=1. In this case $Z((ijklm))=\pi_{g_{lm},g_{kl},g_{jk},g_{ij}}$, and 
the pentagonator is just a 4-cocycle on $G$. As already remarked, this is the 
example given in \cite{Ma99}. In this case 
$\End(I)$ is equivalent to $\Vect$.
\item \#G=1. In this case, without using the funny brackets, we have 
\begin{equation*}
\begin{split}
Z((ijklm))=\ &
(\alpha^1_{h_{ikm},h_{klm},h_{ijk}})
^{-1}\tau_{h_{klm},h_{ijk}}
\alpha^1_{h_{ikm},h_{ijk},h_{klm}}
(\alpha^1_{h_{ijm},h_{jkm},h_{klm}})^{-1}\\
&\times\ \alpha^1_{h_{ijm},h_{jlm},
h_{jkl}}(\alpha^1_{h_{ilm},h_{ijl},
h_{jkl}})^{-1}\alpha^1_{h_{ilm},h_{ikl},h_{ijk}}
\end{split}
\end{equation*}
This is the Crane-Yetter partition function \cite{CY93,CKY97}, which 
they call  the $15j$ symbol, for a finite group instead of a quantum group.  
\end{enumerate} 

We are now ready to define the state-sum, $Z(M,{\cal T})$. Let $v_0$ be 
the number 
of vertices in ${\cal T}$, and $v_1$ the number of edges in ${\cal T}$. In the following 
definition 
the sum is taken over all possible labellings and the product over all 
4-simplices in ${\cal T}$. 
If the orientation of a 4-simplex $S$ induced by the 
ordering 
on its vertices is equal to its orientation induced by the global orientation 
of $M$, then we take $\epsilon(S)=1$. Otherwise we take $\epsilon(S)=-1$. 

\begin{definition}
\label{Z(M,T)}
$$Z(M,{\cal T})=(\#G)^{-v_0}(\#H)^{(v_0-v_1)}\sum_{\ell}\prod_{S}Z(S,\ell)^{\epsilon
(S)}.$$
\end{definition}
Here $Z(S,\ell)$ is defined to be $Z((ijklm))$ for any 4-simplex $S=(ijklm)$ 
in ${\cal T}$. Apart from the extra normalization factor, this is exactly 
the state-sum one obtains from our setup in \cite{Ma99}: since each 
$Z(S,\ell)$ is just one element of $R^*$, rather than a whole matrix of them, 
the tensor product of all these partition functions is just their product, 
$\prod_{S}Z(S,\ell)^{\epsilon(S)}$, which, of course, is just an element of 
$R^*$, so we have $Z(M,{\cal T},\ell)=\prod_{S}Z(S,\ell)^{\epsilon
(S)}$. The quantum dimension of any simple object and any simple 
1-morphism is equal to $1\in R^*$.  

Let us now show that this defines an invariant. 
\begin{theorem}
The state-sum $Z(M,{\cal T})$ is independent of the chosen triangulation 
${\cal T}$.
\end{theorem}
{\bf Proof}. 
We prove invariance under the 4D Pachner moves. 
As explained in the introduction, the two simplicial 4-complexes 
that define a 
4D Pachner move form the boundary of a 5-simplex together. Let us assume 
that this 5-simplex is $(012345)$. By the local semi-flatness conditions, 
the labelling 
of $(012345)$ is uniquely determined by the labels on 
$$(01),\ (12),\ (23),\ 
(34),\ (45)$$
and 
$$(012),\ (013),\ (014),\ (015),\ (023),\ (024),\ (025),\ 
(034),\ (035),\ (045).$$ For short, let us call these labels $g_1,\ldots,\ 
g_5$ and $h_1,\ldots,\ h_{10}$, respectively. 

We first prove the $3\rightleftharpoons 3$ move. The partition function 
corresponding to the left-hand side of this move is equal to
\begin{equation*}
\begin{split} 
Z(&01235)Z(01345)Z(12345)=\lceil(\iota^3_{g_5g_4;g_3;h_1})^{-1}\rceil\\
&\times\lceil\tau_{h_7^{-1}h_9h_5
\alpha^0_{g_5g_4,g_3,g_2g_1};
h_1}\rceil\lceil(\iota^1_{h_7^{-1}h_9h_5\alpha^0_{g_5g_4,g_3,g_2g_1};g_2;g_1})^{-1}\rceil\\
&\times\lceil
\iota^2_{g_5g_4;h_2^{-1}h_5h_1\alpha^0_{g_3,g_2,g_1};g_1}\rceil\lceil\pi_{g_5g_4;g_3;g_2;g_1}\rceil\lceil(\iota^3_{g_5;g_4;h_2})^{-1}\rceil
\lceil\tau_{h_9^{-1}
h_{10}h_8\alpha^0_{g_5,g_4,g_3g_2g_1};
h_2}\rceil\\
&\times\lceil(\iota^1_{h_9^{-1}h_{10}h_8\alpha^0_{g_5,g_4,g_3g_2g_1};g_3g_2;g_1})^{-1}\rceil\lceil
\iota^2_{g_5;h_3^{-1}h_8h_2\alpha^0_{g_4,g_3g_2,g_1};g_1}\rceil
\lceil\pi_{g_5;g_4;g_3g_2;g_1}\rceil\\
&\times
\lceil(\iota^3_{g_5;g_4;h_2^{-1}h_5h_1\alpha^0_{g_3,g_2,g_1}})^{-1}\rceil\lceil\tau_
{h_9^{-1}h_{10}h_8\alpha^0_{g_5,g_4,g_3g_2g_1};
h_2^{-1}h_5h_1\alpha^0_{g_3,g_2,g_1}}\rceil\\
&\times\lceil(\iota^1_{h_9^{-1}h_{10}h_8\alpha^0_{g_5,g_4,g_3g_2g_1};g_3;g_2})^{-1}\rceil
\lceil
\iota^2_{g_5;h_6^{-1}h_8h_5\alpha^0_{g_4,g_3,g_2g_1};g_2}\rceil
\lceil\pi_{g_5;g_4;g_3;g_2}\rceil.
\end{split}
\end{equation*}
On the right-hand side we have 
\begin{equation*}
\begin{split}
Z(&02345)Z(01245)Z(01234)=\lceil(\iota^3_{g_5;g_4;h_5})^{-1}\rceil
\lceil\tau_{h_9^{-1}h_{10}h_8
\alpha^0_{g_5,g_4,g_3g_2g_1};
h_5}\rceil\\
&\times\lceil(\iota^1_{h_9^{-1}h_{10}h_8\alpha^0_{g_5,g_4,g_3g_2g_1};g_3;g_2g_1})^{-1}\rceil\lceil
\iota^2_{g_5;h_6^{-1}h_8h_5\alpha^0_{g_4,g_3,g_2g_1};g_2g_1}\rceil\\
&\times\lceil\pi_{g_5;g_4;g_3;g_2g_1}\rceil
\lceil(\iota^3_{g_5;g_4g_3;h_1})^{-1}\rceil\lceil\tau_{h_7^{-1}h_{10}h_6
\alpha^0_{g_5,g_4g_3,g_2g_1};
h_1}\rceil\\
&\times\lceil(\iota^1_{h_7^{-1}h_{10}h_6\alpha^0_{g_5,g_4g_3,g_2g_1};g_2;g_1})^{-1}\rceil\lceil
\iota^2_{g_5;h_3^{-1}h_6h_1\alpha^0_{g_4g_3,g_2,g_1};g_1}\rceil\lceil\pi_{g_5;g_4g_3;g_2;g_1}\rceil\\
&\times\lceil(\iota^3_{g_4;g_3;h_1})^{-1}\rceil\lceil
\tau_{h_6^{-1}h_8h_5\alpha^0_{g_4,g_3,g_2g_1};
h_1}\rceil\lceil(\iota^1_{h_6^{-1}h_8h_5\alpha^0_{g_4,g_3,g_2g_1};g_2;g_1})^{-1}\rceil\\
&\times\lceil
\iota^2_{g_4;h_2^{-1}h_5h_1\alpha^0_{g_3,g_2,g_1};g_1}\rceil\lceil\pi_{g_4;g_3;g_2;g_1}\rceil.
\end{split}
\end{equation*}
Take the product of the left-hand side with the inverse of the right-hand 
side. After applying $(\rightarrow\otimes\bullet\otimes\bullet\otimes\bullet)$, $(\rightarrow\rightarrow\otimes\bullet\otimes\bullet)$ and the analogous 
identities, and $(\bullet\otimes\bullet\otimes\bullet\otimes\bullet)$ in \ref{sweakmonstruc} we see that this product reduces to
\begin{equation*}
\begin{split} 
&\lceil\tau_{h_7^{-1}h_9h_5\alpha^0_{g_5g_4,g_3,g_2g_1};
h_1}\rceil\lceil\pi_{g_5g_4;g_3;g_2;g_1}\rceil\lceil\tau_{h_9^{-1}
h_{10}h_8\alpha^0_{g_5,g_4,g_3g_2g_1};
h_2}\rceil\lceil\pi_{g_5;g_4;g_3g_2;g_1}\rceil\\
&\times\lceil\tau_
{h_9^{-1}h_{10}h_8\alpha^0_{g_5,g_4,g_3g_2g_1};
h_2^{-1}h_5h_1\alpha^0_{g_3,g_2,g_1}}\rceil
\lceil\pi_{g_5;g_4;g_3;g_2}\rceil 
\lceil\tau^{-1}_{h_9^{-1}h_{10}h_8\alpha^0_{g_5,g_4,g_3g_2g_1};
h_5}\rceil\\
&\times\lceil\pi^{-1}_{g_5;g_4;g_3;g_2g_1}\rceil
\lceil\tau^{-1}_{h_7^{-1}h_{10}h_6
\alpha^0_{g_5,g_4g_3,g_2g_1};
h_1}\rceil
\lceil\pi^{-1}_{g_5;g_4g_3;g_2;g_1}\rceil
\lceil\tau^{-1}_{h_6^{-1}h_8h_5\alpha^0_{g_4,g_3,g_2g_1};h_1}\rceil\\
&\times
\lceil\pi^{-1}_{g_4;g_3;g_2;g_1}\rceil\lceil\iota^2_{g_5;
\alpha^0_{g_4,g_3,g_2};g_1}\rceil
\lceil(\iota^3_{g_5;g_4;
\alpha^0_{g_3,g_2,g_1}})^{-1}\rceil\lceil(\iota^1_{\alpha^0_{g_5,g_4,g_3};
g_2;g_1})^{-1}\rceil\\
&\times
\lceil\tau^{-1}_{h_9^{-1}h_{10}h_8\alpha^0_{g_5,g_4,g_3g_2g_1};
\alpha^0_{g_3,g_2,g_1}}\rceil
\lceil\tau^{-1}_{\alpha^0_{g_5,g_4,g_3};h_1}\rceil.
\end{split}
\end{equation*}
The tensorators, i.e., the $\tau$'s, all cancel because of the relations in  
$(\bullet\otimes\bullet\otimes\bullet)$ and $(\bullet\otimes\bullet\otimes\bullet\otimes\bullet)$. Finally we are left precisely with 
all the terms in relation $(\bullet\otimes\bullet\otimes\bullet\otimes\bullet
\otimes\bullet)$, so we see that our big product is equal to $1$. 

Invariance under the $2\rightleftharpoons 4$ move follows from the 
same calculations. The only difference is that some of the factors on the 
left-hand side now appear at the other side as inverses and vice versa. On 
the right-hand side of the $2\rightleftharpoons 4$ move we have one more 
edge and four more faces than on the left-hand side; in our picture these are 
the edge $(14)$ and the faces $(014),\ (124),\ (134),\ (145)$. Any label of 
$(14)$ is already determined by the the labels of the other edges and the 
local semi-flatness condition. We can choose the label of one of the extra 
faces freely, the labels of the other faces follow from the local semi-flatness 
condition. This means that the product of the factors on the right-hand side 
equals $\#H$ times the product of the factors on the right-hand side. Since we 
normalized our state-sum with the factor $\#G^{-v_0}\#H^{v_0-v_1}$, we get 
the desired result. 

The same kind of argument applies to the $1\rightleftharpoons 5$ move. 
On one side of this move, the one with five arrows, we have one more vertex, 
five more edges, and ten more faces, than on the other side. The labels of one 
of the edges and of four of the faces can be chosen freely. The other 
extra labels are completely determined by local semi-flatness. Again the 
normalization factor ensures invariance.     

We do not prove invariance with respect to the ordering on the 
vertices of $\cal T$ here. As a matter of fact it follows directly from 
the results in \cite{Ma99}, because the proofs of invariance under 
combinatorial isomorphisms do not depend in any way on the assumption that 
$\End(I)$ be $\Vect$. It all follows from the fact that the 2-category 
involved is spherical.  
\qed

\begin{theorem}
Let $\alpha^j,\pi,\tau,\iota^k$ and $(\alpha^j)',\pi',\tau',(\iota^k)'$, for 
$j=1,2$ and $k=1,2,3$, be two 2-equivalent, as defined by Def.~\ref{equival}, 
semi-weak monoidal 2-category structures on ${\mathbb N}(G,H,R)$. Then the 
value of the state-sum using the first 
semi-weak monoidal 2-category structure equals the value of the 
state-sum using the second. 
\end{theorem}
{\bf Proof} We can assume that the automorphisms $G\to G$, $H\to H$ and 
$R\to R$ in Def.~\ref{equival} are all identities, because the state-sum 
is taken over all labellings. Let $\mu,\Phi,\phi,\psi,\chi$ define the 
2-equivalence. The 
labellings of the two state-sums correspond in the following way 
\begin{equation*}
\begin{split}
g_{ij}'&\leftrightarrow g_{ij}\\
h'_{ijk}&\leftrightarrow h_{ijk}\Phi_{g_{jk},g_{ij}}
\end{split}
\end{equation*}
Note that under this correspondence we have 
\begin{equation*}
\begin{split}
h'_{jkl}{h'}^{-1}_{ikl}{h'}_{ijl}{h'}^{-1}_{ijk}&=
\Phi^{-1}_{g_{jk},g_{ij}}\Phi^{-1}_{g_{kl},g_{jk}g_{ij}}
h_{jkl}h^{-1}_{ikl}h_{ijl}h^{-1}_{ijk}
\Phi_{g_{kl}g_{jk},g_{ij}}\Phi_{g_{kl},g_{jk}}\\
&=\Phi^{-1}_{g_{jk},g_{ij}}\Phi^{-1}_{g_{kl},g_{jk}g_{ij}}
\alpha^0_{g_{kl},g_{jk},g_{ij}}
\Phi_{g_{kl}g_{jk},g_{ij}}\Phi_{g_{kl},g_{jk}}\\
&=(\alpha^0)'_{g'_{kl},g'_{jk},g'_{ij}}.
\end{split}
\end{equation*}
Let us now have a look at the partition function $Z(ijklm)'$:
\begin{equation*}
\begin{split}
Z(ijklm)'=&\lceil ((\iota^3)'_{g'_{lm},g'_{kl},h'_{ijk}})^{-1}\rceil\lceil
(\tau')_{h'_{klm},h'_{ijk}}\rceil\lceil((\iota^1)'_
{h'_{klm},g'_{jk},g'_{ij}})^{-1}\rceil\\
&\times\lceil(\iota^2)'_{g'_{lm},h'_{jkl},g'_{ij}}\rceil
\lceil\pi'_{g'_{lm},g'_{kl},g'_{jk},g'_{ij}}\rceil.
\end{split}
\end{equation*}
Using the correspondence between the labellings, which we defined above, and 
the identities in Def.~\ref{equival}, $Z(ijklm)'$ is seen to be equal to 
\begin{equation*}
\begin{split}
&Z(ijklm)\lceil \chi^{-1}_{g_{lm}g_{kl},h_{ijk}} 
\rceil\lceil \chi_{g_{kl},h_{ijk}}\rceil\lceil 
\chi_{g_{lm},h_{ijk}\alpha^0_{g_{kl},g_{jk},g_{ij}}}
\rceil\\
&\times\lceil \psi^{-1}_{h_{klm},g_{jk}}\rceil\lceil 
\psi^{-1}_{h_{klm}(\alpha^0)^{-1}_{g_{lm},g_{kl},
g_{jk}},g_{ij}}
\rceil\lceil \psi_{h_{klm},g_{jk}g_{ij}}\rceil\lceil
\phi_{g_{lm}g_{kl},g_{jk},g_{ij}}\rceil\\
&\times\lceil
\phi_{g_{lm},g_{kl},g_{jk}g_{ij}}\rceil\lceil
\phi^{-1}_{g_{lm},g_{kl},g_{jk}}\rceil\lceil
\phi^{-1}_{g_{lm},g_{kl}g_{jk},g_{ij}}\rceil\lceil
\phi^{-1}_{g_{kl},g_{jk},g_{ij}}\rceil.
\end{split}
\end{equation*}
By local semi-flatness this equals
\begin{equation*}
\begin{split}
&Z(ijklm)\lceil \chi^{-1}_{g_{lm}g_{kl},h_{ijk}} 
\rceil\lceil \chi_{g_{kl},h_{ijk}}\rceil\lceil 
\chi_{g_{lm},h_{ijl}}\rceil\lceil\chi^{-1}_{g_{lm},h_{ikl}}\rceil\lceil
\chi_{g_{lm},h_{jkl}}
\rceil\\
&\times\lceil \psi^{-1}_{h_{klm},g_{jk}}\rceil\lceil 
\psi^{-1}{h_{jlm},g_{ij}}\rceil\lceil\psi_{h_{jkm},g_{ij}}\rceil\lceil
\psi^{-1}_{h_{jkl},g_{ij}}
\rceil\lceil \psi_{h_{klm},g_{jk}g_{ij}}\rceil\lceil
\phi_{g_{lm}g_{kl},g_{jk},g_{ij}}\rceil\\
&\times\lceil
\phi_{g_{lm},g_{kl},g_{jk}g_{ij}}\rceil\lceil
\phi^{-1}_{g_{lm},g_{kl},g_{jk}}\rceil\lceil
\phi^{-1}_{g_{lm},g_{kl}g_{jk},g_{ij}}\rceil\lceil
\phi^{-1}_{g_{kl},g_{jk},g_{ij}}\rceil.
\end{split}
\end{equation*}
Note that every factor other than $Z(ijklm)$ in the formula above corresponds 
to a tetrahedron in ${\cal T}$, rather than a 4-simplex. 
It is well known that in a 4-manifold without boundary each tetrahedron 
belongs to the boundary of exactly two 4-simplices, appearing with a positive 
sign in one boundary and with a negative sign in the other. The reason is 
that the link of each vertex is homeomorphic to $S^4$~\cite{RS82}. Therefore 
all factors other than 
$Z(ijklm)$ cancel out, because in the product over all 4-simplices in the 
state-sum they appear exactly twice, once with a positive and once with 
a negative sign. 
\qed
\vskip.2cm 

One nice consequence of our 
approach via 2-categories is that our construction generalizes several known 
constructions at once. If we take the trivial monoidal 2-category structure 
on ${\mathbb N}(G,H,R)$, we recover Yetter's~\cite{Ye93} invariants 
corresponding to homotopy 2-types. Porter~\cite{Po99} generalized Yetter's 
construction using homotopy $n$-types, for arbitrary $n$. However, his 
state-sums for homotopy 3-types are different from ours. For a given 
triangulated 4-manifold, $M=(M,{\cal T})$, and a given homotopy 3-type, 
Porter's construction yields a state-sum which simply counts the number of 
possible labellings of $\cal T$ up to some normalization. Note that 
these labellings are not equal to ours, because Porter also assigns labels to 
the tetrahedra. 
   
For $H=\{1\}$, and $R={\mathbb C}$, our partition function 
is defined by a 4-cocycle on $G$. Birmingham and Rakowski~\cite{BR96} show 
that for $G={\mathbb Z}/n{\mathbb Z}$, with $n$ a non-negative integer, the 
invariant is 
equal to Yetter's~\cite{Ye93} untwisted invariant, because the product of the 
4-cocycles is always equal to $1$ for a closed 4-manifold. 

We already mentioned that for $\# G=1$ we get the 
Crane-Yetter~\cite{CY93,CKY97} invariants for finite groups. This case has 
been studied by Birmingham and Rakowski in \cite{BR94} for 
$H={\mathbb Z}/n{\mathbb Z}$, for $n$ a non-negative integer, 
and $R={\mathbb C}$. The model 
that they study corresponds 
to the case in which only $\tau$ in Definition~\ref{sweakmonstruc} is 
non-trivial. They show that the 
partition function in their case can be obtained by evaluation of the 
intersection form defined on the second cohomology group of the simplicial 
complex $\cal T$ that defines the triangulation, with coefficients in 
${\mathbb Z}/n{\mathbb Z}$, against the fundamental homology cycle of $\cal T$. In our 
context their definition of $\tau$ becomes:
$$\tau_{h_1,h_2}=\exp({\frac{2\pi ik}{n}}[h_1h_2]).$$
Here $0 < k < n-1$ is an integer and $[h_1h_2]$ is defined to be 
$h_1h_2\bmod n$. The $h_i$ can be defined as the integers $0,\ldots,n-1$. 
Birmingham and Rakowski also 
present explicit calculations of the state-sum for the complex projective 
plane, ${\mathbb C}P^2$, for $n=2,3$. The values they obtained are:
$$Z(\pm{\mathbb C}P^2)=0\ \  \mbox{for $n=2$},$$
$$Z(\pm{\mathbb C}P^2)=\mp 3\sqrt{3}i\ \ \mbox{for $n=3$}.$$ This shows that 
the invariant is non-trivial.

In \cite{BR95} Birmingham and Rakowski present a construction of 4-manifold 
invariants that correspond to ours for the case in which 
$G=H={\mathbb Z}/n{\mathbb Z}$, with $n$ a non-negative integer, and $R={\mathbb C}$, 
and only 
$\iota^1$ is 
taken to be non-trivial in Definition~\ref{sweakmonstruc}. In our context 
their definition of $\iota^1$ becomes:
$$\iota^1_{h,g_1,g_2}=\exp({\frac{2\pi ik}{n^2}} h(g_1+g_2-[g_1+g_2])).$$ 
Here $0 < k < n^2$ is an integer and $[g_1+g_2]$ is defined to be 
$g_1+g_2\bmod n$. Also in this definition we take $g_i$ and $h_i$ to be the 
integers $0,\ldots,n-1$. In \cite{BR95} Birmingham and Rakowski 
calculate the state-sums for ${\mathbb R}P^3\times S^1$, $S^4$, $S^3\times S^1$, 
and $L(5,1)$, a lens space. We recall the values they obtained:
$$Z({\mathbb R}P^3\times S^1)=\left\{ \begin{array}{ll}
2\cdot 2^{\delta_2(k)} & \mbox{for $n$ even} \\
1 &\mbox{otherwise}\end{array}\right., $$
$$Z(S^4)=1,$$ 
$$Z(S^3\times S^1)=1,$$
$$Z(L(5,1))=\left\{ \begin{array}{ll}
5\cdot 5^{\delta_5(k)} & \mbox{for $n\equiv 0\pmod 5$} \\ 
1 & \mbox{otherwise}\end{array}\right..$$
The $\pmod n$ delta function, $\delta_n$, is defined by  
$$\delta_n(k)=\left\{ \begin{array}{ll}
1 & \mbox{if $k\equiv 0\pmod n$} \\
0 & \mbox{otherwise}\end{array}\right..$$
These computations show that the invariants are rather non-trivial. Birmingham and Rakowski~\cite{BR94} mention that 
it would be interesting to do similar computations for the case in which 
one multiplies the above mentioned partition functions, i.e., $\iota^1$ and 
$\delta$. Here we have set everything in a more general context, 
thereby providing one point of view for all the different models 
that Birmingham and Rakowski consider. In our partition function we also have 
a factor $\iota^2$. Looking at Birmingham and Rakowski's examples 
it is not hard to find an example of $\iota^2$ in the same context. We 
can define 
$$\iota^2_{g_1,h,g_3}=\exp({\frac{2\pi ik}{n}}[g_1hg_3]).$$
One could take the product of $\tau$, $\iota^1$, and $\iota^2$, for the 
partition function, as a special case of our construction.    
         
\section{Postnikov systems}
\label{P}
The connection between equivalence classes of semi-weak monoidal 2-category structures on 
${\mathbb N}(G,H,R)$ and Postnikov 
systems, as sketched in this section, is based on 
the conjecture that a semi-weak monoidal 2-category, as 
defined in Def.~\ref{sweakmonstruc}, can be seen as a semi-weak 3-category, 
as defined by Tamsamani~\cite{Ta95}, with one object. 
 
As remarked in the introduction already, several 
people~\cite{BD98,Bat97,Ta95} have suggested a definition of weak 
$n$-categories. Unfortunately 
the question whether these definitions are ``equivalent'' is extremely subtle 
and has 
not been settled yet. Tamsamani follows an approach via simplicial sets 
which stays 
very close to the ideas coming from homotopy theory. Since we want to relate 
semi-weak monoidal 2-category structures to Postnikov systems, Tamsamani's setup is 
convenient here.  
Tamsamani shows that his definition of a category correponds to the 
``ordinary'' definition. He also 
shows that his definition of a weak 2-category is equivalent to the definition of a 
bi-category as defined in \cite{Be67}, which is the definition that underlies Gordon, Power and 
Street's definition of a tricategory. It is therefore very reasonable to 
conjecture that 
a weak 3-category in the sense of Tamsamani's definition yields a tricategory and vice versa. 
However, the verification of this conjecture would take many pages, as can be seen from the 
length of Tamsamani's proof of the equivalence of the definitions of 
his weak 2-categories and Benabou's~\cite{Be67} bi-categories. Therefore we 
do not attempt to prove the conjecture here. We mean this section 
to be motivational for the earlier parts of this chapter and are, for that 
reason, also a little sketchy in this section.  

All definitions of weak $n$-categories are complicated and inductive, so we do not 
wish to repeat Tamsamani's definition here. As a matter of fact we only need a consequence of 
his results, which we explain now. In the second part of his PhD dissertation 
\cite{Ta95} Tamsamani 
realizes an idea that was first sketched by Grothendieck~\cite{Gr83}. 
Tamsamani defines weak $n$-groupoids for any $n\in {\mathbb N}$, which are weak $n$-categories of which all $k$-morphisms 
are invertible up to higher order equivalences, and 
shows that equivalence classes of weak $n$-groupoids correspond bijectively to 
homotopy classes of $n$-anticonnected CW-complexes. An $n$-anticonnected space is one for which all 
homotopy groups of order greater than $n$ vanish. Here equivalence is again a very 
subtle matter. Under this correspondence equivalence classes of $k$-morphisms, 
with $0\leq k\leq n$, correspond exactly to the elements of the $k$th-order homotopy group. 
Our definition of ${\mathbb N}(G,H,R)$ is just the 'linearized' version of 
a 3-groupoid with one object of which the 1-morphisms are the elements of 
$G$, the 2-morphisms the elements of $H$, and the 3-morphisms the elements of 
$R^*$. One could call such a 3-groupoid a {\it groupal 2-groupoid}. In our 
case the actions of $G$ on $H$ and $R^*$ are trivial. 
Therefore, the 
equivalence classes of the structures 
of a semi-weak monoidal 2-category, i.e., a weak 3-category with one object in the sense 
of Tamsamani's definition with strict units and ``trivial'' tensor product 
of simple objects with simple 1- and 2-morphisms, on ${\mathbb N}(G,H,R)$ 
correspond bijectively to homotopy classes of CW-complexes of which the only 
non-vanishing homotopy groups are $\pi_1=G$, $\pi_2=H$, and $\pi_3=R^*$, 
and for which the actions of $\pi_1$ on $\pi_2$ and $\pi_3$ are trivial. Such 
CW-complexes we call connected 3-anticonnected $>1$-simple. This is 
analogous to the results stated in \cite{Qu98}.  The proof in Quinn's paper of 
the analogous result for monoidal groupoids is essentially due to 
\cite{BFSV98,Fie98,FQ93}. In this text we 
put more emphasis on the connection with higher dimensional algebra.

Thus the classification of semi-weak monoidal 
2-category structures 
on ${\mathbb N}(G,H,R)$ boils down, conjecturally, to the classification up 
to homotopy equivalence of connected 3-anticonnected $>1$-simple 
CW-complexes $E$ with 
$\pi_1(E)=G$, $\pi_2(E)=H$, $\pi_3(E)=R^*$. 
It is well known~\cite{Bre93,Wh78} that such a classification is obtained by 
the theory of Postnikov systems. Some people may not be familiar with 
this theory, so let us briefly 
sketch its key idea. Let $X$ be a connected $n$-anticonnected $>1$-simple 
CW-complex with $n\geq 1$, and let $A$ be an abelian group. Then there is a 
one-to-one correspondence 
between the homotopy equivalence classes of connected $(n+1)$-anticonnected 
$>1$-simple 
CW-complexes $Y$ 
of which all homotopy groups up to order $n$ coincide with those of $X$ and of which 
$\pi_{n+1}(Y)$ is equal to $A$, and homotopy classes of maps $\alpha\colon X\to K(A,n+2)$. 
Here $K(A,n+2)$ is the so called Eilenberg-MacLane space of order $n+2$ with group A, 
of which the only non-vanishing homotopy group is $\pi_{n+2}=A$. If $A$ is finite, 
then $K(A,n+2)$ is equal to $B_A^{n+2}$. Here we define $B_A$ to be the classifying 
space of $A$, 
which is a simplicial group itself, and we 
define inductively $B_A^{n}=B_{B_A^{n-1}}$. Given such a 
map $\alpha$, one can take $Y$ to be a CW-approximation of the principal 
fibration induced by $\alpha$, which is the pull back, over $\alpha$, of the 
so called path-loop fibration of $K(A,n+2)$. As a set this principal 
fibration is given by $\{(x,\gamma)\in X\times P(K(A,n+2))\vert \alpha(x)=\gamma(1)\}$, where $P(K(A,n+2))$ is the space of all paths in $K(A,n+2)$ which 
start at the base-point. Note that the fibre of the principal fibration is 
$\Omega(K(A,n+2))$, the space of all loops starting and ending at the 
base-point in $K(A,n+2)$, which, as is well known, is weakly homotopy 
equivalent to $K(A,n+1)$. Conversely, one can prove 
that any $Y$ of the aforementioned type is a CW-approximation of the 
principal fibration induced by such a map. Thus two maps $\alpha\colon B_G\to B_H^3$ and $\beta\colon 
W(\alpha)\to B_{R^*}^4$, where $W(\alpha)$ is a CW-approximation of the 
principal fibration induced by $\alpha$, correspond to   
a $>1$-simple connected 3-anticonnected CW-complex $E$ with $\pi_1(E)=G$, $\pi_2(E)=H$, and 
$\pi_3(E)=R^*$, that is unique up to homotopy equivalence. Since homotopy classes of maps 
$X\to K(A,n+2)$ correspond bijectively to cohomology classes in 
$H^{n+2}(X,A)$ (see \cite{Bre93,Wh78}), we 
arrive at the following conjecture:

\begin{conjecture}
\label{conj}
The equivalence classes of semi-weak monoidal 2-category structures on 
${\mathbb N}(G,H,R)$, as defined in Defs.~\ref{sweakmonstruc} and 
\ref{equival}, correspond bijectively 
to pairs of cohomology classes $\alpha\in H^3(B_G,H)$ and 
$\beta=\beta_{\alpha}\in H^4(W(\alpha),R^*)$. 
\end{conjecture}
One implication of the conjecture is easy to prove 
directly: 
our partition function on one 4-simplex in Def.~\ref{Z((ijklm))} defines 
exactly a 4-cocycle on $W(\alpha)$ with values in $R^*$. A simplicial set $S=S(\alpha)$ 
whose geometric realization yields $W(\alpha)$ is due to Blakers and was 
worked out by Brown and Higgins, see \cite{Bro99} and 
references therein, or, 
equivalently, by the geometric realization of the nerve of the groupal 
groupoid corresponding to $\alpha$, see \cite{Ta95} and references therein. 
In this construction, there is only one 0-simplex in $S$, a 1-simplex for 
every $g\in G$, 
a 2-simplex for every triple $g_1,g_2,g_3\in G$ satisfying $g_1g_2=g_3$ 
and every $h\in H$, a 3-simplex for every quadruple of 2-simplices in $S$ such 
that the edges match up around the 3-simplex and such that the four elements 
$h_1,h_2,h_3,h_4\in H$ corresponding to the triangles satisfy  
$$h_1h_2^{-1}h_3h_4^{-1}=\alpha^0_{g_3,g_2,g_1}.$$
Note that this description corresponds exactly to the diagram we have drawn in Fig.~\ref{label}. 
In general there is an $n$-simplex for every $n+1$-tuple of $n-1$-simplices 
in $S$ whose faces match up appropriately. A 4-cocycle on $W(\alpha)$ with 
values in $R^*$ is nothing but a map from the Abelian group generated by the 
4-cells in the cellular model of $W(\alpha)$ to $R^*$, the instances of which 
multiply up to one around the boundary of a 5-simplex. Our partition function 
in Def.~\ref{Z((ijklm))} is an example of such a map and its invariance 
under the $3\rightleftharpoons 3$ Pachner move shows that it satisfies the 
4-cocycle condition. 

It would be very nice if we could derive the structural maps in 
Definition~\ref{sweakmonstruc} 
directly from the Postnikov invariants, thereby proving the conjecture. 
Unfortunately we have not been able to do this completely.  
It is clear that $\alpha$ represents the 0-associator, $\alpha^0$. It is also easy to obtain 
the 1-associator, $\alpha^1$, and the tensorator, $\tau$. Let $i\colon B_H^2\to W(\alpha)$ be 
the embedding of the fibre in the fibration, then $i^*\colon H^4(W(\alpha),R^*)\to H^4(B_H^2,R^*)$ 
defines the pull-back $i^*(\beta)\in H^4(B_H^2,R^*)$. Quinn~\cite{Qu98} 
showed that the cohomology 
classes in $H^4(B_H^2,R^*)$ correspond bijectively to the equivalence classes of weak braided monoidal 
structures on ${\mathbb N}(H,R)$. 
In our context these are exactly the weak braided monoidal structures on $\End(1)$, the category of endomorphisms of 
the identity object. Quinn also shows explicitly how a 4-cocycle representing 
an element in $H^4(B_H^2,R^*)$ consists of a 3-cocycle and a 
2-cochain on $H$ 
with values in $R^*$ which satisfy the hexagon equations in Definition~\ref{sweakmonstruc}. The 
3-cocycle represents the 1-associator, $\alpha^1$, and the 2-cochain represents the tensorator, 
$\tau$. The present author does not know how to obtain the remaining maps and 
relations in Definition~\ref{sweakmonstruc} from an element in 
$H^4(B_H^2,R^*)$.

For any $\beta\in H^4(W(\alpha),R^*)$ we can define invariants of 
4-manifolds in the same way as 
Quinn does in \cite{Qu95}. Given $f\colon M\to W(\alpha)$ one can evaluate the pull-back 
$f^*(\beta)\in H^4(M,R^*)$ on the fundamental homology class of the 
4-manifold, $[M]$. For any $\alpha'$ cohomologous to $\alpha$, any $\beta'$ 
cohomologous to $\beta$, and any $g\colon M\to W(\alpha')$ homotopic to 
$f$, we have $g^*(\beta')([M])=f^*(\beta)([M])$. Quinn then takes a certain 
weighted sum of $f^*(\beta)([M])$ over an arbitrary set of 
representatives of all 
homotopy classes of maps $f$. For a precise definition of the weights see 
\cite{Qu95}. As a matter of fact, Quinn only works out concretely his very 
abstract 
construction, which he defines in any dimension, for the classifying 
space of a finite group and cocycles on 
that space. In dimension four that corresponds to the restricted 
case in which $H=\{1\}$. The interesting invariants in 
dimension 4, that we have sketched above, were never considered by Quinn, 
or anyone else, explicitly. It is clear, by the arguments following our 
conjecture, that our state-sum invariants for a given semi-weak monoidal 
structure are equal to Quinn's invariants for the corresponding Postnikov 
invariants. This shows 
immediately 
that our state-sum invariants are homotopy invariants, rather than PL 
invariants. This is not surprising given the fact that we use finite groups. 
In dimension 3, the Dijkgraaf-Witten~\cite{DW90} invariants are homotopy 
invariants as well. As already mentioned, the Turaev-Viro invariants are 
real homeomorphism invariants, but they require the use of quantum groups 
instead of finite groups. The categorical construction that underlies 
the Dijkgraaf-Witten and the Turaev-Viro invariants is the same though; it 
is the specific input in that construction that makes the difference. 
It is therefore reasonable to look for categorifications of the quantum 
groups, the representations of which could be the right input in our 
construction, presented in \cite{Ma99}, for obtaining true PL invariants. This idea led 
Crane and Frenkel~\cite{CF94} to the definition of a Hopf category in the 
first place.

\section{Final remarks}

First of all let us address the question of examples. We already mentioned
at the end of Section \ref{Z} that Birmingham and 
Rakowski's~\cite{BR94,BR95,BR96} constructions can be seen as special cases of 
our construction. Therefore, their computations show that there are 
non-trivial examples of the kind of invariant that we describe in this 
paper. It remains to be seen if there are more examples. Section~\ref{P} 
indicates that there should be many more examples, since any homotopy 3-type 
theoretically provides one example. In \cite{Bro91} Brown has computed some 
homotopy 3-types using non-abelian tensor products of groups. Hopefully his 
results will provide more concrete examples of our construction. 
   
We can ask ourselves how powerful we can expect our state-sum invariants 
to be. 
By the 'conjectural' relation with Postnikov systems and the relation with 
Freed and Quinn's work 
in \cite{FQ93,Qu95} it is clear that our invariants are homotopy invariants rather than PL invariants. Depending on one's point of view one can find this 
disappointing or encouraging. We take the latter point of view, because 
the TQFT programme, as sketched in \cite{BD95}, for example, remains still to 
be developed in dimension four. Any interesting examples of four-dimensional 
TQFT's, even of a homotopic nature, are welcome at this stage of 
the development of the TQFT programme. For the case $G=1$ and $H={\bf Z}/n{\bf Z}$ Birmingham 
and Rakowski~\cite{BR94} have shown that the partition function can be 
obtained by the evaluation of the intersection form defined on the second 
cohomology group of the triangulation against the fundamenal form of the 
manifold. It would be interesting to know if there are any relations 
between our invariants and other classical homotopy invariants. 

As already mentioned in the introduction, there is another construction of 
4-manifold invariants using finite groups: the 
Crane-Frenkel~\cite{CF94} construction for the categorification of the quantum 
double of a finite group~\cite{CY97}. This has been worked out in detail by Carter, Kauffman, 
and Saito in~\cite{CKS99}. It would be worthwile to figure out the precise 
relation between that construction and ours. In~\cite{Ma99} we conjectured that 
the 2-category of representations of an involutory Hopf category is a 
spherical 2-category, and that the Crane-Frenkel construction using involutory 
Hopf categories and our 
construction using the 2-categories of representations of Hopf categories yield 
the same invariants. However, as mentioned in the introduction of this paper, 
we did assume that $\End(I)$ is Vect in \cite{Ma99}. How the two constructions 
relate to one another when $\End(I)$ is an arbitrary tortile category we do 
not know. This is certainly something to be investigated and a good point to 
start would be the case involving only finite groups. 
   
The final remark we want to make is that there should be results, that are analogous to 
our results in this paper, for braided monoidal 2-categories. 
The general definition of these 2-categorical structures was first given by Kapranov and 
Voevodsky~\cite{KV94}. Later Baez and Neuchl~\cite{BN96} and Crans~\cite{Cr98} corrected some flaws in that definition. 
In \cite{BD95} Baez and Dolan conjectured that braided monoidal 2-categories are 4-categories with 
one object and one 1-morphism. Let us assume that this is true for a moment. In that case we 
see, from Tamsamani's~\cite{Ta95,Ta96} results, that braided monoidal structures on 
$\N(G,H,R)$ correspond to connected CW-complexes 
of which $\pi_2=G,\ \pi_3=H,\ \pi_4=R^*$ are the only non-vanishing homotopy groups. By the 
theory of Postnikov systems we see that these CW-complexes are classified up to homotopy 
equivalence by two cohomology classes, $\alpha\in H^4(K(G,2),H)$ and $\beta\in H^5(W(\alpha),R^*)$, where $W(\alpha)$ is the path-loop fibration over 
$K(G,2)$ induced by $\alpha$.
Note the shift in the order of the cohomology groups. It would be nice to work out concretely 
all the maps and relations that define braided monoidal structures on $\N(G,H,R)$, 
analogously to what we do in this paper, and write down the invariants that Baez and Langford 
\cite{BL982,BL98} defined. Also in that case it would be desirable to find arguments by which one can extract all these maps and relations directly from 
the cohomology classes.

\section{Acknowledgements}

First of all I thank my supervisor, Prof. L. Crane, for the many 
discussions we have had about the results that I present in this paper and 
his encouragement. I also thank Prof. D. Yetter and Prof. F. Quinn 
for helpful e-mail discussions. Prof. J. Baez I thank for his information 
about results in the theory of $n$-categories that are relevant to my work 
and Prof. T. Porter I thank for his information 
about relevant results in homotopy theory.  

\bibliography{TQFT.bib}

\end{document}